\documentclass[oneside,english,a4paper]{amsart}
\usepackage[T1]{fontenc}
\usepackage[latin9]{inputenc}
\usepackage{amstext}
\usepackage{amsthm}
\usepackage{amssymb}

\makeatletter
\numberwithin{equation}{section}
\numberwithin{figure}{section}
\theoremstyle{plain}
\newtheorem{thm}{\protect\theoremname}
\theoremstyle{remark}
\newtheorem{rem}[thm]{\protect\remarkname}
\theoremstyle{plain}
\newtheorem{lem}[thm]{\protect\lemmaname}
\theoremstyle{definition}
\newtheorem{defn}[thm]{\protect\definitionname}
\theoremstyle{plain}
\newtheorem{prop}[thm]{\protect\propositionname}
\theoremstyle{remark}
\newtheorem*{rem*}{\protect\remarkname}

\usepackage{hyperref}   

\makeatother

\usepackage{babel}
\providecommand{\definitionname}{Definition}
\providecommand{\lemmaname}{Lemma}
\providecommand{\propositionname}{Proposition}
\providecommand{\remarkname}{Remark}
\providecommand{\theoremname}{Theorem}

\begin{document}
\title[doubly perturbed Yamabe problems]{Compactness and blow up results for doubly perturbed Yamabe problems
on manifolds with umbilic boundary}
\author{Marco G. Ghimenti}
\address{M. G. Ghimenti, \newline Dipartimento di Matematica Universit\`a di Pisa
Largo B. Pontecorvo 5, 56126 Pisa, Italy}
\email{marco.ghimenti@unipi.it}

\author{Anna Maria Micheletti}
\address{A. M. Micheletti, \newline Dipartimento di Matematica Universit\`a di Pisa
Largo B. Pontecorvo 5, 56126 Pisa, Italy}
\email{a.micheletti@dma.unipi.it.}

\begin{abstract}
Given a compact Riemannian manifold with umbilic boundary, the Yamabe
boundary problem studies if there exist conformal scalar-flat metrics
such that $\partial M$ has constant mean curvature. In this paper
we address to the stability of this problem with respect of perturbation
of mean curvature of the boundary and scalar curvature of the manifold.
In particular we prove that the Yamabe boundary problem is stable
under perturbation of the mean curvature and the scalar curvature
from below, while it is not stable if one of the two curvatures is
perturbed from above. 
\end{abstract}

\maketitle

\section{Introduction}

Let $(M,g)$, a smooth, compact Riemannian manifold of dimension $n\ge3$
with boundary. In \cite{Es} Escobar asked it there exists a conformal
metric $\tilde{g}=u^{\frac{4}{n-2}}g$ for which $M$ has zero scalar
curvature and constant boundary mean curvature. 

This problem can be understood as a generalization of the Riemann
mapping theorem and it is equivalent to finding a positive solution
to the following nonlinear boundary value problem

\begin{equation}
\left\{ \begin{array}{cc}
L_{g}u=0 & \text{ in }M\\
B_{g}u+(n-2)u^{\frac{n}{n-2}}=0 & \text{ on }\partial M
\end{array}\right..\label{eq:Pconf}
\end{equation}
Where $L_{g}=\Delta_{g}-\frac{n-2}{4(n-1)}R_{g}$ and $B_{g}=-\frac{\partial}{\partial\nu}-\frac{n-2}{2}h_{g}$
are respectively the conformal Laplacian and the conformal boundary
operator, $R_{g}$ is the scalar curvature of the manifold, $h_{g}$
is the mean curvature of the $\partial M$ and $\nu$ is the outer
normal with respect to $\partial M$ . If $M$ is of \emph{positive
type}, that is when 
\[
Q(M,\partial M):=\inf_{u\in H^{1}\smallsetminus0}\frac{\int\limits _{M}\left(|\nabla u|^{2}+\frac{n-2}{4(n-1)}R_{g}u^{2}\right)dv_{g}+\int\limits _{\partial M}\frac{n-2}{2}h_{g}u^{2}d\sigma_{g}}{\left(\int\limits _{\partial M}|u|^{\frac{2(n-1)}{n-2}}d\sigma_{g}\right)^{\frac{n-2}{n-1}}}
\]
is strictly positive, equation (\ref{eq:Pconf}) could have multiple
solutions, and the question of compactness of solution arises naturally.
In fact, if the boundary of $M$ is umbilic, and the Weyl tensor $W_{g}$
never vanishes on the boundary, the full set of solution of (\ref{eq:Pconf})
is compact. This is proved in \cite{GM20}, for dimensions $n>8$,
and in \cite{GMsub}, for dimensions $n=6,7,8$. We recall that the
boundary of $M$ is called umbilic if the trace-free second fundamental
form of $\partial M$ is zero everywhere.

Also, the authors show in \cite{GMdcds} that the problem is stable
for perturbation from below of the mean curvature, while in \cite{GMP}
with Pistoia they prove that there is a blow up phenomenon when perturbing
the mean curvature from above. This recalls a similar result from
the Yamabe problems on boundaryless manifolds, in which perturbations
form below of the scalar curvature preserve the compactness of the
set of solutions (see \cite{dru,DH}). 

At this point it is interesting to study what happens when one perturbs
both the scalar and the mean curvature, and to investigate compactness
versus blow up of solutions in this framework. Thus, we study the
linearly perturbed problem
\begin{equation}
\left\{ \begin{array}{cc}
-\Delta_{g}u+\frac{n-2}{4(n-1)}R_{g}u+\varepsilon_{1}\alpha u=0 & \text{ in }M\\
\frac{\partial u}{\partial\nu}+\frac{n-2}{2}h_{g}u+\varepsilon_{2}\beta u=(n-2)u^{\frac{n}{n-2}} & \text{ on }\partial M
\end{array}\right.\label{eq:Prob-2}
\end{equation}
or, in a more compact form, 
\[
\left\{ \begin{array}{cc}
L_{g}u-\varepsilon_{1}\alpha u=0 & \text{ in }M\\
B_{g}u-\varepsilon_{2}\beta u+(n-2)u^{\frac{n}{n-2}}=0 & \text{ on }\partial M
\end{array}\right.,
\]
 where $\varepsilon_{1},\varepsilon_{2}$ are small positive parameters
and $\alpha,\beta:M\rightarrow\mathbb{R}$ are smooth functions. Here
we choose $\varepsilon_{1}$ sufficiently such that $-L_{g}+\varepsilon_{1}\alpha$
is still a positive definite operator. 

Our aim is to prove that, if we linearly perturb the mean curvature
term $h_{g}$ with a negative smooth function, and jointly we perturb
the scalar curvature term $R_{g}$ with another negative smooth function,
the set of solution is still compact. On the contrary, if one between
scalar and mean curvature is perturbed from above, the compactness
of solutions is lost. Our main results read as
\begin{thm}
\label{thm:main}Let $(M,g)$ a smooth, $n$-dimensional Riemannian
manifold of positive type not conformally equivalent to the standard
ball with regular umbilic boundary $\partial M$. 

Let $\alpha,\beta:M\rightarrow\mathbb{R}$ smooth functions such that
$\alpha,\beta<0$ on $\partial M$. Suppose that $n\ge8$ and that
the Weyl tensor $W_{g}$ is not vanishing on $\partial M$. Then,
there exists a positive constant $C$ such that for any $\varepsilon_{1},\varepsilon_{2}\ge0$
small enough and for any $u>0$ solution of (\ref{eq:Prob-2}) it
holds 
\[
C^{-1}\le u\le C\text{ and }\|u\|_{C^{2,\eta}(M)}\le C
\]
for some $0<\eta<1$. The constant $C$ does not depend on $u,\varepsilon_{1},\varepsilon_{2}$. 
\end{thm}

\begin{thm}
\label{thm:main2}Let $(M,g)$ a smooth, $n$-dimensional Riemannian
manifold of positive type not conformally equivalent to the standard
ball with regular umbilic boundary $\partial M$. 

Let $\alpha,\beta:M\rightarrow\mathbb{R}$ smooth functions. Suppose
that $n\ge8$ and that the Weyl tensor $W_{g}$ is not vanishing on
$\partial M$. If $\alpha>0$ on $\partial M$ or $\beta>0$ on $\partial M$,
then there exists a sequence of solutions $u_{\varepsilon_{1},\varepsilon_{2}}$
of (\ref{eq:Prob-2}) which blows up at a point of the boundary when
$(\varepsilon_{1},\varepsilon_{2})\rightarrow(0,0)$.
\end{thm}

Let us shortly comment these two results.
\begin{itemize}
\item In a series of paper \cite{dru,DH,DHR} Druet, Hebey and Robert studied
the stability of classical Yamabe problem under perturbation of scalar
curvature terms. They proved that the set of solutions of $-\Delta_{g}u+\frac{n-2}{4(n-1)}a(x)u=cu^{\frac{n+2}{n-2}}\text{ in }M$
is compact if $a(x)\le R_{g}(x)$ on $M$, thus the problem is stable
perturbing $R_{g}$ from below, while they found counterexamples to
compactness when $a(x)$ is greater than $R_{g}(x)$. In \cite{GMdcds,GMP}
the same problem is studied in the case of boundary Yamabe equations
by perturbing the mean curvature term and a matching compactness versus
blow up phenomenon appears. So there is a strong analogy between the
role of $R_{g}$ in classical case and $h_{g}$ in boundary case.
We continue here the same analysis, by perturbing both the curvature
terms at the same time, to complete the study. It appears that the
problem is stable only when perturbing both terms with non positive
functions, while it is enough to perturb from above one between $h_{g}$
and $R_{g}$ to lose compactness of the solutions.
\item We worked here in the framework of umbilic boundary manifolds. In
a recent paper \cite{GM-nonumb}, we studied the case of manifold
with non umbilic boundary, that is when the trace-free second fundamental
form is non zero in any point of $\partial M$. In this case it is
possible to have compactness also for positive small perturbation
of the scalar curvature. We want to remind that in the case of non
umbilic boundary the compactness of solution for the unperturbed case
was proved by Almaraz \cite{Al} and Kim, Musso and Wei \cite{KMW}. 
\item In the unperturbed case the compactness of solutions for umbilic manifolds
has been proved for dimensions $n\ge6$ (see \cite{GM20,GMsub}).
It should be possible to apply the same technique also in the perturbed
case to extend Theorem \ref{thm:main} to $n=6,7$. It is less clear
to us if Theorem \ref{thm:main2} could be extended to lower dimensions.
Case $n=5$ remains open also for the unperturbed problem.
\item Our theorems consider only perturbations that are everywhere positive
or everywhere negative on $M$. However, in Remark \ref{rem:sign-changing}
it is shown that it is possible to construct a sign changing $\alpha$
such that for any $\beta$ the set of solutions in no more compact,
or a sign changing $\beta$ such that for any $\alpha$ the set of
solutions in no more compact. We do not know if it would be possible
to craft a sign changing perturbation for which compactness is preserved.
\end{itemize}
The paper is organized as follows. Hereafter we recall some basic
definitions and all the preliminary notions useful to achieve the
result. Section \ref{sec:The-compactness-result} is devoted to the
proof of the compactness theorem, while in Section \ref{sec:The-non-compactness}
we prove the non compactness result.

\subsection{Notations and preliminary definitions}
\begin{rem}[Notations]
We will use the indices $1\le i,j,k,m,p,r,s\le n-1$ and $1\le a,b,c,d\le n$.
Moreover we use the Einstein convention on repeated indices. We denote
by $g$ the Riemannian metric, by $R_{abcd}$ the full Riemannian
curvature tensor, by $R_{ab}$ the Ricci tensor and by $R_{g}$ and
$h_{g}$ respectively the scalar curvature of $(M,g)$ and the mean
curvature of $\partial M$; moreover the Weyl tensor of $(M,g)$ will
be denoted by $W_{g}$. The bar over an object (e.g. $\bar{W}_{g}$)
will means the restriction to this object to the metric of $\partial M$.

Finally, on the half space $\mathbb{R}_{+}^{n}=\left\{ y=(y_{1},\dots,y_{n-1},y_{n})\in\mathbb{R}^{n},\!\ y_{n}\ge0\right\} $
we set $B_{r}(y_{0})=\left\{ y\in\mathbb{R}^{n},\!\ |y-y_{0}|\le r\right\} $
and $B_{r}^{+}(y_{0})=B_{r}(y_{0})\cap\left\{ y_{n}>0\right\} $.
When $y_{0}=0$ we will use simply $B_{r}=B_{r}(y_{0})$ and $B_{r}^{+}=B_{r}^{+}(y_{0})$.
On the half ball $B_{r}^{+}$ we set $\partial'B_{r}^{+}=B_{r}^{+}\cap\partial\mathbb{R}_{+}^{n}=B_{r}^{+}\cap\left\{ y_{n}=0\right\} $
and $\partial^{+}B_{r}^{+}=\partial B_{r}^{+}\cap\left\{ y_{n}>0\right\} $.
On $\mathbb{R}_{+}^{n}$ we will use the following decomposition of
coordinates: $(y_{1},\dots,y_{n-1},y_{n})=(\bar{y},y_{n})=(z,t)$
where $\bar{y},z\in\mathbb{R}^{n-1}$ and $y_{n},t\ge0$.

Fixed a point $q\in\partial M$, we denote by $\psi_{q}:B_{r}^{+}\rightarrow M$
the Fermi coordinates centered at $q$. We denote by $B_{g}^{+}(q,r)$
the image of $B_{r}^{+}$. When no ambiguity is possible, we will
denote $B_{g}^{+}(q,r)$ simply by $B_{r}^{+}$, omitting the chart
$\psi_{q}$.
\end{rem}

We introduce the following notation for integral quantities which
recur often in the paper
\[
I_{m}^{\alpha}:=\int_{0}^{\infty}\frac{s^{\alpha}ds}{\left(1+s^{2}\right)^{m}}.
\]
By direct computation (see \cite[Lemma 9.4]{Al}) it holds
\begin{align}
I_{m}^{\alpha}=\frac{2m}{\alpha+1}I_{m+1}^{\alpha+2} & \text{ for }\alpha+1<2m\label{eq:Iam}\\
I_{m}^{\alpha}=\frac{2m}{2m-\alpha-1}I_{m+1}^{\alpha} & \text{ for }\alpha+1<2m\nonumber \\
I_{m}^{\alpha}=\frac{2m-\alpha-3}{\alpha+1}I_{m}^{\alpha+2} & \text{ for }\alpha+3<2m\nonumber 
\end{align}
Also, we have the following integral identities:
\begin{align}
\int_{0}^{\infty}\frac{t^{k}dt}{(1+t)^{m}} & =\frac{k!}{(m-1)(m-2)\cdots(m-1-k)}\label{eq:t-integrali}\\
\int_{0}^{\infty}\frac{dt}{(1+t)^{m}} & =\frac{1}{m-1}\nonumber 
\end{align}
and, by change of variables
\begin{equation}
\int_{\mathbb{R}_{+}^{n}}\frac{|\bar{y}|^{\alpha}y_{n}^{\beta}}{\left[(1+y_{n})^{2}+|\bar{y}|^{2}\right]^{\gamma}}d\bar{y}dy_{n}=\omega_{n-2}I_{\gamma}^{\alpha+n-2}\int_{0}^{\infty}\frac{y_{n}^{\beta}}{(1+y_{n})^{2\gamma-\alpha-n+1}}dy_{n}\label{eq:doppioint}
\end{equation}
where $\omega_{n-2}$ is the volume of $\mathbb{S}^{n-1}$.

We shortly recall here the well known function ${\displaystyle U(y):=\frac{1}{\left[(1+y_{n})^{2}+|\bar{y}|^{2}\right]^{\frac{n-2}{2}}}}$
which is also called the standard bubble and which is the unique solution,
up to translations and rescaling, of the nonlinear critical problem
. 
\begin{equation}
\left\{ \begin{array}{ccc}
-\Delta U=0 &  & \text{on }\mathbb{R}_{+}^{n};\\
\frac{\partial U}{\partial y_{n}}=-(n-2)U^{\frac{n}{n-2}} &  & \text{on \ensuremath{\partial}}\mathbb{R}_{+}^{n}.
\end{array}\right.\label{eq:Udelta}
\end{equation}
We set 
\begin{equation}
j_{l}:=\partial_{l}U=-(n-2)\frac{y_{l}}{\left[(1+y_{n})^{2}+|\bar{y}|^{2}\right]^{\frac{n}{2}}}\label{eq:jl}
\end{equation}
\[
\partial_{k}\partial_{l}U=(n-2)\left\{ \frac{ny_{l}y_{k}}{\left[(1+y_{n})^{2}+|\bar{y}|^{2}\right]^{\frac{n+2}{2}}}-\frac{\delta^{kl}}{\left[(1+y_{n})^{2}+|\bar{y}|^{2}\right]^{\frac{n}{2}}}\right\} 
\]
\begin{equation}
j_{n}:=y^{b}\partial_{b}U+\frac{n-2}{2}U=-\frac{n-2}{2}\frac{|y|^{2}-1}{\left[(1+y_{n})^{2}+|\bar{y}|^{2}\right]^{\frac{n}{2}}}.\label{eq:jn}
\end{equation}
and we recall that $j_{1},\dots,j_{n}$ are a base of the space of
the $H^{1}$ solutions of the linearized problem 
\begin{equation}
\left\{ \begin{array}{ccc}
 & -\Delta\phi=0 & \text{on }\mathbb{R}_{+}^{n},\\
 & \frac{\partial\phi}{\partial t}+nU^{\frac{2}{n-2}}\phi=0 & \text{on \ensuremath{\partial}}\mathbb{R}_{+}^{n},\\
 & \phi\in H^{1}(\mathbb{R}_{+}^{n}).
\end{array}\right.\label{eq:linearizzato}
\end{equation}
Given a point $q\in\partial M$, we introduce now the function $\gamma_{q}$
which arises from the second order term of the expansion of the metric
$g$ on $M$ (see \ref{eq:gij}). The choice of this function plays
a twofold role in this paper. On the one hand, using the function
$\gamma_{q}$ we are able to perform the estimates of Lemmas \ref{lem:coreLemma},
\ref{lem:taui} and Proposition \ref{prop:stimawi}. On the other
hand, it gives the correct correction to the standard bubble in order
to perform finite dimensional reduction.

For the proof of the following Lemma we refer to \cite[Lemma 3]{GMP}
and \cite[Proposition 5.1]{Al}.
\begin{lem}
\label{lem:vq}Assume $n\ge5$. Given a point $q\in\partial M$, there
exists a unique $\gamma_{q}:\mathbb{R}_{+}^{n}\rightarrow\mathbb{R}$
a solution of the linear problem 
\begin{equation}
\left\{ \begin{array}{ccc}
-\Delta\gamma=\left[\frac{1}{3}\bar{R}_{ijkl}(q)y_{k}y_{l}+R_{ninj}(q)y_{n}^{2}\right]\partial_{ij}^{2}U &  & \text{on }\mathbb{R}_{+}^{n}\\
\frac{\partial\gamma}{\partial y_{n}}=-nU^{\frac{2}{n-2}}\gamma &  & \text{on }\partial\mathbb{R}_{+}^{n}
\end{array}\right.\label{eq:vqdef}
\end{equation}
which is $L^{2}(\mathbb{R}_{+}^{n})$-orthogonal to the functions
$j_{1},\dots,j_{n}$ defined in (\ref{eq:jl}) and (\ref{eq:jn}).

Moreover it holds
\begin{equation}
|\nabla^{\tau}\gamma_{q}(y)|\le C(1+|y|)^{4-\tau-n}\text{ for }\tau=0,1,2.\label{eq:gradvq}
\end{equation}
\begin{equation}
\int_{\mathbb{R}_{+}^{n}}\gamma_{q}\Delta\gamma_{q}dy\le0,\label{new}
\end{equation}

\begin{equation}
\int_{\partial\mathbb{R}_{+}^{n}}U^{\frac{n}{n-2}}(t,z)\gamma_{q}(t,z)dz=0\label{eq:Uvq}
\end{equation}
\begin{equation}
\gamma_{q}(0)=\frac{\partial\gamma_{q}}{\partial y_{1}}(0)=\dots=\frac{\partial\gamma_{q}}{\partial y_{n-1}}(0)=0.\label{eq:dervq}
\end{equation}

Finally the map $q\mapsto\gamma_{q}$ is $C^{2}(\partial M)$.
\end{lem}

\subsection{Expansion of the metric\label{sec:Expansion}}

Since the boundary $\partial M$ is umbilic, given $q\in\partial M$
there exists a conformally related metric $\tilde{g}_{q}=\Lambda_{q}^{\frac{4}{n-2}}g$
such that some geometric quantities at $q$ have a simpler form which
will be summarized later in this paragraph. We have 
\[
\Lambda_{q}(q)=1,\ \frac{\partial\Lambda_{q}}{\partial y_{k}}(q)=0\text{ for all }k=1,\dots,n-1.
\]
 Also, we have that $\tilde{u}_{q}:=\Lambda_{q}u$, is a solution
of (\ref{eq:Prob-2}) if and only if $u$ solves the following problem
\begin{equation}
\left\{ \begin{array}{cc}
L_{\tilde{g}_{q}}u-\varepsilon_{1}\tilde{\alpha}u=0 & \text{ in }M\\
B_{\tilde{g}_{q}}u+(n-2)u^{\frac{n}{n-2}}-\varepsilon_{2}\tilde{\beta}u=0 & \text{ on }\partial M
\end{array}\right..\label{eq:P-conf}
\end{equation}
where  $\tilde{\alpha}:=\Lambda_{q}^{-\frac{4}{n-2}}\alpha$ and $\tilde{\beta}:=\Lambda_{q}^{-\frac{2}{n-2}}\beta$. 

In the following expansion and in section \ref{sec:The-compactness-result},
in order to simplify notations, we will omit the \emph{tilde} symbols,
since we will always work in the conformal metric $\tilde{g}$, while
in Section \ref{sec:The-non-compactness} we will switch between metrics,
so we will keep $g$ and $\tilde{g}$ explicitly indicated. 

With this metric we have the following expansions.
\begin{rem}
\label{rem:confnorm}In Fermi conformal coordinates around $q\in\partial M$,
it holds (see \cite{M1})
\begin{equation}
|\text{det}g_{q}(y)|=1+O(|y|^{N})\text{ with }N\text{ arbitrarily large}\label{eq:|g|}
\end{equation}
\begin{eqnarray}
|h_{ij}(y)|=O(|y^{4}|) &  & |h_{g}(y)|=O(|y^{4}|)\label{eq:hij}
\end{eqnarray}
\begin{align}
g_{q}^{ij}(y)= & \delta^{ij}+\frac{1}{3}\bar{R}_{ikjl}y_{k}y_{l}+R_{ninj}y_{n}^{2}\label{eq:gij}\\
 & +\frac{1}{6}\bar{R}_{ikjl,m}y_{k}y_{l}y_{m}+R_{ninj,k}y_{n}^{2}y_{k}+\frac{1}{3}R_{ninj,n}y_{n}^{3}\nonumber \\
 & +\left(\frac{1}{20}\bar{R}_{ikjl,mp}+\frac{1}{15}\bar{R}_{iksl}\bar{R}_{jmsp}\right)y_{k}y_{l}y_{m}y_{p}\nonumber \\
 & +\left(\frac{1}{2}R_{ninj,kl}+\frac{1}{3}\text{Sym}_{ij}(\bar{R}_{iksl}R_{nsnj})\right)y_{n}^{2}y_{k}y_{l}\nonumber \\
 & +\frac{1}{3}R_{ninj,nk}y_{n}^{3}y_{k}+\frac{1}{12}\left(R_{ninj,nn}+8R_{nins}R_{nsnj}\right)y_{n}^{4}+O(|y|^{5})\nonumber 
\end{align}
\begin{equation}
R_{g_{q}}(y)=O(|y|^{2})\text{ and }\partial_{ii}^{2}R_{g_{q}}=-\frac{1}{6}|\bar{W}|^{2}\label{eq:Rii}
\end{equation}
\begin{equation}
\partial_{tt}^{2}\bar{R}_{g_{q}}=-2R_{ninj}^{2}-2R_{ninj,ij}\label{eq:Rtt}
\end{equation}
\begin{equation}
\bar{R}_{kl}=R_{nn}=R_{nk}=R_{nn,kk}=0\label{eq:Ricci}
\end{equation}
\begin{equation}
R_{nn,nn}=-2R_{nins}^{2}.\label{eq:Rnnnn}
\end{equation}
All the quantities above are calculate in $q\in\partial M$, unless
otherwise specified.
\end{rem}

If we choose $\varepsilon_{1}$ sufficiently small in order to have
that $-L_{g}+\varepsilon_{1}\alpha$ is a positive definite operator,
we can define an equivalent scalar product on $H^{1}$ as 
\begin{equation}
\left\langle \left\langle u,v\right\rangle \right\rangle _{g}=\int_{M}\left(\nabla_{g}u\nabla_{g}v+\frac{n-2}{4(n-1)}R_{g}uv+\varepsilon_{1}\alpha uv\right)d\mu_{g}+\frac{n-2}{2}\int_{\partial M}h_{g}uvd\nu_{g}\label{eq:prodscal}
\end{equation}
which leads to the norm $\|\cdot\|_{g}$ equivalent to the usual one. 

With this norm we have that $\Lambda_{q}$ is an isometry. In fact,
by (\ref{eq:prodscal}), for any $u,v\in H^{1}(M)$,
\[
\left\langle \left\langle \Lambda_{q}u,\Lambda_{q}v\right\rangle \right\rangle _{g}=\left\langle \left\langle u,v\right\rangle \right\rangle _{\tilde{g}_{q}}\text{ and, consequently, }\|\Lambda_{q}u\|_{g}=\|u\|_{\tilde{g}_{q}}.
\]

\subsection{Variational framework.}

Given $1\le t\le\frac{2(n-1)}{n-2}$ we have the well known embedding
\[
i:H^{1}(M)\rightarrow L^{t}(\partial M).
\]
We define, by the scalar product $\left\langle \left\langle \cdot,\cdot\right\rangle \right\rangle _{g}$,
\[
i_{\alpha}^{*}:L^{t}(\partial M)\rightarrow H^{1}(M)
\]
in the following sense: given $f\in L^{\frac{2(n-1)}{n-2}}(\partial M)$
there exists a unique $v\in H^{1}(M)$ such that 
\begin{align}
v=i_{\alpha}^{*}(f) & \iff\left\langle \left\langle v,\varphi\right\rangle \right\rangle _{g}=\int_{\partial M}f\varphi d\sigma\text{ for all }\varphi\label{eq:istella}\\
 & \iff\left\{ \begin{array}{ccc}
-\Delta_{g}v+\frac{n-2}{4(n-1)}R_{g}v+\varepsilon_{1}\alpha=0 &  & \text{on }M;\\
\frac{\partial v}{\partial\nu}+\frac{n-2}{2}h_{g}v=f &  & \text{on \ensuremath{\partial}}M.
\end{array}\right.\nonumber 
\end{align}
So Problem (\ref{eq:Prob-2}) is equivalent to find $v\in H^{1}(M)$
such that 

\[
v=i_{\alpha}^{*}(f(v)-\varepsilon_{2}\beta v)
\]
where 
\[
f(v)=(n-2)\left(v^{+}\right)^{\frac{n}{n-2}}.
\]
Notice that, if $v\in H_{g}^{1}$, then $f(v)\in L^{\frac{2(n-1)}{n}}(\partial M)$. 

Problem (\ref{eq:Prob-2}) has also a variational structure and a
positive solution for (\ref{eq:Prob-2}) is a critical point for the
following functional defined on $H^{1}(M)$ 
\begin{align*}
J_{\varepsilon_{1},\varepsilon_{2},g}(v)=J_{g}(v): & =\frac{1}{2}\int_{M}|\nabla_{g}v|^{2}+\frac{n-2}{4(n-1)}R_{g}v^{2}+\varepsilon_{1}\alpha v^{2}d\mu_{g}+\frac{n-2}{4}\int_{\partial M}h_{g}v^{2}d\sigma_{g}\\
 & +\frac{1}{2}\int_{\partial M}\varepsilon_{2}\beta v^{2}d\sigma_{g}-\frac{(n-2)^{2}}{2(n-1)}\int_{\partial M}\left(v^{+}\right)^{\frac{2(n-1)}{n-2}}d\sigma_{g}.
\end{align*}
We remark that, defined
\begin{align*}
\tilde{J}_{\varepsilon_{1},\varepsilon_{2},\tilde{g}}(v)=\tilde{J}_{\tilde{g}}(v): & =\frac{1}{2}\int_{M}|\nabla_{\tilde{g}}v|^{2}+\frac{n-2}{4(n-1)}R_{\tilde{g}}v^{2}+\varepsilon_{1}\tilde{\alpha}v^{2}d\mu_{\tilde{g}}+\frac{n-2}{4}\int_{\partial M}h_{\tilde{g}}v^{2}d\sigma_{\tilde{g}}\\
 & +\frac{1}{2}\int_{\partial M}\varepsilon_{2}\tilde{\beta}v^{2}d\sigma_{\tilde{g}}-\frac{(n-2)^{2}}{2(n-1)}\int_{\partial M}\left(v^{+}\right)^{\frac{2(n-1)}{n-2}}d\sigma_{\tilde{g}},
\end{align*}
we have 
\begin{equation}
J_{g}(\Lambda_{q}u)=\tilde{J}_{\tilde{g}}(u).\label{eq:Jlambda}
\end{equation}

Given $q\in\partial M$ and $\psi_{q}^{\partial}:\mathbb{R}_{+}^{n}\rightarrow M$
the Fermi coordinates in a neighborhood of $q$; we define 
\begin{align*}
W_{\delta,q}(\xi) & =U_{\delta}\left(\left(\psi_{q}^{\partial}\right)^{-1}(\xi)\right)\chi\left(\left(\psi_{q}^{\partial}\right)^{-1}(\xi)\right)=\\
 & =\frac{1}{\delta^{\frac{n-2}{2}}}U\left(\frac{y}{\delta}\right)\chi(y)=\frac{1}{\delta^{\frac{n-2}{2}}}U\left(x\right)\chi(\delta x)
\end{align*}
where $y=(z,t)$, with $z\in\mathbb{R}^{n-1}$ and $t\ge0$, $\delta x=y=\left(\psi_{q}^{\partial}\right)^{-1}(\xi)$
and $\chi$ is a radial cut off function, with support in ball of
radius $R$. In an analogous way, given $\gamma_{q}$ as in Lemma
\ref{lem:vq} we define 
\[
V_{\delta,q}(\xi)=\frac{1}{\delta^{\frac{n-2}{2}}}\gamma_{q}\left(\frac{1}{\delta}\left(\psi_{q}^{\partial}\right)^{-1}(\xi)\right)\chi\left(\left(\psi_{q}^{\partial}\right)^{-1}(\xi)\right)
\]
 and, given $j_{a}$ defined in (\ref{eq:jl}) and (\ref{eq:jn})
we define 
\[
Z_{\delta,q}^{b}(\xi)=\frac{1}{\delta^{\frac{n-2}{2}}}j_{b}\left(\frac{1}{\delta}\left(\psi_{q}^{\partial}\right)^{-1}(\xi)\right)\chi\left(\left(\psi_{q}^{\partial}\right)^{-1}(\xi)\right).
\]
By means of $\left\langle \left\langle \cdot,\cdot\right\rangle \right\rangle _{g}$
it is possible to decompose $H^{1}$ in the direct sum of the following
two subspaces 
\begin{align*}
\tilde{K}_{\delta,q} & =\text{Span}\left\langle \Lambda_{q}Z_{\delta,q}^{1},\dots,\Lambda_{q}Z_{\delta,q}^{n}\right\rangle \\
\tilde{K}_{\delta,q}^{\bot} & =\left\{ \varphi\in H^{1}(M)\ :\ \left\langle \left\langle \varphi,\Lambda_{q}Z_{\delta,q}^{b}\right\rangle \right\rangle _{g}=0,\ b=1,\dots,n\right\} 
\end{align*}
and to define the projections 
\[
\tilde{\Pi}=H^{1}(M)\rightarrow\tilde{K}_{\delta,q}\text{ and }\tilde{\Pi}^{\bot}=H^{1}(M)\rightarrow\tilde{K}_{\delta,q}^{\bot}.
\]
 Notice that, since $\Lambda_{q}$ is an isometry, we have that $\varphi\in\tilde{K}_{\delta,q}$
if and only if $\Lambda_{q}^{-1}\varphi\in K_{\delta,q}$ and the
same holds for $\tilde{K}_{\delta,q}^{\bot}$. 

In Section \ref{sec:The-non-compactness}, we will look for a solution
$\tilde{u}_{q}=\Lambda_{q}u$ of (\ref{eq:Prob-2}) which has the
form 
\[
\tilde{u}_{q}=\Lambda_{q}\left(W_{\delta,q}+\delta^{2}V_{\delta,q}+\phi\right)=\tilde{W}_{\delta,q}+\delta^{2}\tilde{V}_{\delta,q}+\tilde{\phi}
\]
where $\tilde{\phi}\in\tilde{K}_{\delta,q}^{\bot}$. By means of $i_{\alpha}^{*}$
this is equivalent to the following pair of equations 
\begin{align}
\tilde{\Pi}\left\{ \tilde{W}_{\delta,q}+\delta^{2}\tilde{V}_{\delta,q}+\tilde{\phi}-i_{\alpha}^{*}\left[f(\tilde{W}_{\delta,q}+\delta^{2}\tilde{V}_{\delta,q}+\tilde{\phi})-\varepsilon_{2}\beta(\tilde{W}_{\delta,q}+\delta^{2}\tilde{V}_{\delta,q}+\tilde{\phi})\right]\right\}  & =0\label{eq:Pi}\\
\tilde{\Pi}^{\bot}\left\{ \tilde{W}_{\delta,q}+\delta^{2}\tilde{V}_{\delta,q}+\tilde{\phi}-i_{\alpha}^{*}\left[f(\tilde{W}_{\delta,q}+\delta^{2}\tilde{V}_{\delta,q}+\tilde{\phi})-\varepsilon_{2}\beta(\tilde{W}_{\delta,q}+\delta^{2}\tilde{V}_{\delta,q}+\tilde{\phi})\right]\right\}  & =0.\label{eq:Pibot}
\end{align}

\section{The compactness result\label{sec:The-compactness-result}}

In this section, firstly we recall a Pohozaev type identity which
will gives us a fundamental sign condition to rule out the possibility
of blowing up sequence (see subsection \ref{sec:Sign-estimates}).
A recall of preliminary results on blow up points is collected in
subsection \ref{sec:Isolated-and-simple}, while a careful analysis
of blow up sequences is performed in subsection \ref{sec:Blowup-estimates}.
The proof of Theorem \ref{thm:main} is given in subsection \ref{subsec:The-finite-dimensional}.
Throughout this section we work in $\tilde{g}$ metric. For the sake
of readability we will omit the \emph{tilde} symbol throughout this
section.

\subsection{A Pohozaev type identity\label{sec:Pohozaev}}

We will use this version of a local Pohozaev type identity (see \cite{Al,GM20}).
\begin{thm}[Pohozaev Identity]
\label{thm:poho} Let $u$ a $C^{2}$-solution of the following problem
\[
\left\{ \begin{array}{cc}
-\Delta_{g}u+\frac{n-2}{4(n-1)}R_{g}u+\varepsilon_{1}\alpha=0=0 & \text{ in }B_{r}^{+}\\
\frac{\partial u}{\partial\nu}+\frac{n-2}{2}h_{g}u+\varepsilon_{2}\beta u=(n-2)u^{\frac{n}{n-2}} & \text{ on }\partial'B_{r}^{+}
\end{array}\right.
\]
for $B_{r}^{+}=\psi_{q}^{-1}(B_{g}^{+}(q,r))$ for $q\in\partial M$.
Let us define
\[
P(u,r):=\int\limits _{\partial^{+}B_{r}^{+}}\left(\frac{n-2}{2}u\frac{\partial u}{\partial r}-\frac{r}{2}|\nabla u|^{2}+r\left|\frac{\partial u}{\partial r}\right|^{2}\right)d\sigma_{r}+\frac{r(n-2)^{2}}{2(n-1)}\int\limits _{\partial(\partial'B_{r}^{+})}u^{\frac{2(n-1)}{n-2}}d\bar{\sigma}_{g},
\]
and
\begin{multline*}
\hat{P}(u,r):=-\int\limits _{B_{r}^{+}}\left(y^{a}\partial_{a}u+\frac{n-2}{2}u\right)[(L_{g}-\Delta)u]dy+\frac{n-2}{2}\int\limits _{\partial'B_{r}^{+}}\left(\bar{y}^{k}\partial_{k}u+\frac{n-2}{2}u\right)h_{g}ud\bar{y}\\
+\varepsilon_{1}\int\limits _{B_{r}^{+}}\left(y^{a}\partial_{a}u+\frac{n-2}{2}u\right)\alpha udy\\
+\frac{n-2}{2}\varepsilon_{2}\int\limits _{\partial'B_{r}^{+}}\left(\bar{y}^{k}\partial_{k}u+\frac{n-2}{2}u\right)\beta ud\bar{y}.
\end{multline*}
Then $P(u,r)=\hat{P}(u,r)$.

Here $a=1,\dots,n$, $k=1,\dots,n-1$ and $y=(\bar{y},y_{n})$, where
$\bar{y}\in\mathbb{R}^{n-1}$ and $y_{n}\ge0$.
\end{thm}

\subsection{Isolated and isolated simple blow up points\label{sec:Isolated-and-simple}}

Here we recall the definitions of some type of blow up points, and
we give the basic properties about the behavior of these blow up points
(see \cite{Al,FA,HL,M3}). We will omit the proofs of the well known
results.

Let $\left\{ u_{i}\right\} _{i}$ be a sequence of positive solution
to 
\begin{equation}
\left\{ \begin{array}{cc}
L_{g_{i}}u-\varepsilon_{1,i}\alpha_{i}u=0 & \text{ in }M\\
B_{g_{i}}u+(n-2)u^{\frac{n}{n-2}}-\varepsilon_{2,i}\beta_{i}u=0 & \text{ on }\partial M
\end{array}\right..\label{eq:Prob-i}
\end{equation}
where $\alpha_{i}=\Lambda_{x_{i}}^{-\frac{4}{n-2}}\alpha\rightarrow\Lambda_{x_{0}}^{-\frac{4}{n-2}}\alpha$,
$\beta_{i}=\Lambda_{x_{i}}^{-\frac{2}{n-2}}\beta\rightarrow\Lambda_{x_{0}}^{-\frac{2}{n-2}}\beta$,
$x_{i}\rightarrow x_{0}$, $g_{i}\rightarrow g_{0}$ in the $C_{\text{loc}}^{3}$
topology and $0\le\varepsilon_{1,i},\varepsilon_{2,i}<\bar{\varepsilon}$.
\begin{defn}
\label{def:blowup}

1) We say that $x_{0}\in\partial M$ is a blow up point for the sequence
$u_{i}$ of solutions of (\ref{eq:Prob-i}) if there is a sequence
$x_{i}\in\partial M$ of local maxima of $\left.u_{i}\right|_{\partial M}$
such that $x_{i}\rightarrow x_{0}$ and $u_{i}(x_{i})\rightarrow+\infty.$

Shortly we say that $x_{i}\rightarrow x_{0}$ is a blow up point for
$\left\{ u_{i}\right\} _{i}$. 

2) We say that $x_{i}\rightarrow x_{0}$ is an \emph{isolated} blow
up point for $\left\{ u_{i}\right\} _{i}$ if $x_{i}\rightarrow x_{0}$
is a blow up point for $\left\{ u_{i}\right\} _{i}$ and there exist
two constants $\rho,C>0$ such that
\[
u_{i}(x)\le Cd_{\bar{g}}(x,x_{i})^{\frac{2-n}{2}}\text{ for all }x\in\partial M\smallsetminus\left\{ x_{i}\right\} ,\ d_{\bar{g}}(x,x_{i})<\rho.
\]

Given $x_{i}\rightarrow x_{0}$ an isolated blow up point for $\left\{ u_{i}\right\} _{i}$,
and given $\psi_{i}:B_{\rho}^{+}(0)\rightarrow M$ the Fermi coordinates
centered at $x_{i}$, we define the spherical average of $u_{i}$
as
\[
\bar{u}_{i}(r)=\frac{2}{\omega_{n-1}r^{n-1}}\int_{\partial^{+}B_{r}^{+}}u_{i}\circ\psi_{i}d\sigma_{r}
\]
and
\[
w_{i}(r):=r^{\frac{2-n}{2}}\bar{u}_{i}(r)
\]
for $0<r<\rho.$

3) We say that $x_{i}\rightarrow x_{0}$ is an \emph{isolated simple}
blow up point for $\left\{ u_{i}\right\} _{i}$ solutions of (\ref{eq:Prob-i})
if $x_{i}\rightarrow x_{0}$ is an isolated blow up point for $\left\{ u_{i}\right\} _{i}$
and there exists $\rho$ such that $w_{i}$ has exactly one critical
point in the interval $(0,\rho)$.
\end{defn}

Given $x_{i}\rightarrow x_{0}$ a blow up point for $\left\{ u_{i}\right\} _{i}$,
we set
\[
M_{i}:=u_{i}(x_{i})\ \text{ and }\ \delta_{i}:=M_{i}^{\frac{2}{2-n}}.
\]
Obviously $M_{i}\rightarrow+\infty$ and $\delta_{i}\rightarrow0$.

We recall two propositions whose proofs can be found in \cite{A3}
and in \cite{FA}.
\begin{prop}
\label{prop:4.1}Let $x_{i}\rightarrow x_{0}$ is an \emph{isolated}
blow up point for $\left\{ u_{i}\right\} _{i}$ and $\rho$ as in
Definition \ref{def:blowup}. We set 
\[
v_{i}(y)=M_{i}^{-1}(u_{i}\circ\psi_{i})(M_{i}^{\frac{2}{2-n}}y),\text{ for }y\in B_{\rho M_{i}^{\frac{n-2}{2}}}^{+}(0).
\]
Then, given $R_{i}\rightarrow\infty$ and $c_{i}\rightarrow0$, up
to subsequences, we have
\begin{enumerate}
\item $|v_{i}-U|_{C^{2}\left(B_{R_{i}}^{+}(0)\right)}<c_{i}$;
\item ${\displaystyle \lim_{i\rightarrow\infty}\frac{R_{i}}{\log M_{i}}=0}$.
\end{enumerate}
\end{prop}

\begin{prop}
\label{prop:Lemma 4.4}Let $x_{i}\rightarrow x_{0}$ be an isolated
simple blow-up point for $\left\{ u_{i}\right\} _{i}$. Let $\eta$
small. If $0<\bar{\varepsilon}\le1$ is small enough and $0\le\varepsilon_{1,i},\varepsilon_{2,i}\le\bar{\varepsilon}$,
there exist $C,\rho>0$ such that 
\[
M_{i}^{\lambda_{i}}|\nabla^{k}u_{i}(\psi_{i}(y))|\le C|y|^{2-k-n+\eta}
\]
for $y\in B_{\rho}^{+}(0)\smallsetminus\left\{ 0\right\} $ and $k=0,1,2$.
Here $\lambda_{i}=\left(\frac{2}{n-2}\right)(n-2-\eta)-1$.
\end{prop}

\begin{prop}
\label{prop:eps-i}Let $x_{i}\rightarrow x_{0}$ be an isolated simple
blow-up point for $\left\{ u_{i}\right\} _{i}$ and $\alpha,\beta<0$.
Then $\varepsilon_{1,i}\delta_{i}\rightarrow0$ and $\varepsilon_{2,i}\rightarrow0$.
\end{prop}

\begin{proof}
We compute the Pohozaev identity in a ball of radius $r$ and we set
$\frac{r}{\delta_{i}}=:R_{i}\rightarrow\infty$. We estimate any term
of $P(u_{i},r_{i})$ and $\hat{P}(u_{i},r_{i})$.

Set 
\begin{align*}
I_{1}(u,r):= & \int\limits _{\partial^{+}B_{r}^{+}}\left(\frac{n-2}{2}u\frac{\partial u}{\partial r}-\frac{r}{2}|\nabla u|^{2}+r\left|\frac{\partial u}{\partial r}\right|^{2}\right)d\sigma_{r}\\
I_{2}(u,r):= & \frac{r(n-2)^{2}}{2(n-1)}\int\limits _{\partial(\partial'B_{r}^{+})}u^{\frac{2(n-1)}{n-2}}d\bar{\sigma}_{g},
\end{align*}
we have $P(u_{i},r)=I_{1}(u_{i},r)+I_{2}(u_{i},r)$

By Proposition \ref{prop:Lemma 4.4} we obtain
\begin{align*}
I_{1}(u_{i},r) & =M_{i}^{-2\lambda_{i}}I_{1}(M_{i}^{\lambda_{i}}u_{i},r)\le cM_{i}^{-2\lambda_{i}}\int\limits _{\partial^{+}B_{r}^{+}}|y|^{2(2-n+\eta)}d\sigma_{r}\le c\delta_{i}^{\lambda_{i}(n-2)};\\
I_{2}(u_{i},r) & \le cM^{-\lambda_{i}\frac{2(n-1)}{n-2}}\le c\delta_{i}^{\lambda_{i}(n-2)}.
\end{align*}
Then

\begin{equation}
P(u_{i},r)\le\delta_{i}^{\lambda_{i}(n-2)}.\label{eq:poho1}
\end{equation}
In a similar way we decompose
\begin{align*}
\hat{P}(u_{i},r): & =-\int\limits _{B_{r}^{+}}\left(y^{a}\partial_{a}u+\frac{n-2}{2}u\right)[(L_{g}-\Delta)u]dy+\frac{n-2}{2}\int\limits _{\partial'B_{r}^{+}}\left(\bar{y}^{k}\partial_{k}u+\frac{n-2}{2}u\right)h_{g}ud\bar{y}\\
 & +\varepsilon_{1}\int\limits _{B_{r}^{+}}\left(y^{a}\partial_{a}u+\frac{n-2}{2}u\right)\alpha udy\\
 & +\frac{n-2}{2}\varepsilon_{2}\int\limits _{\partial'B_{r}^{+}}\left(\bar{y}^{k}\partial_{k}u+\frac{n-2}{2}u\right)\beta ud\bar{y}=:I_{3}(u,r)+I_{4}(u,r)+I_{5}(u,r)+I_{6}(u,r).
\end{align*}
The terms $I_{3},I_{4}$ and $I_{6}$ are been estimated in \cite{GMdcds}.
For the sake of completeness we report here the main steps of the
estimates. By Proposition \ref{prop:Lemma 4.4} and by definition
of $v_{i}$ we have 
\[
|\nabla^{k}v_{i}(s)|\le M_{i}^{\eta\frac{2}{n-2}}|1+|s||^{2-k-n}=\delta_{i}^{-\eta}|1+|s||^{2-k-n}.
\]
So, recalling that $|h_{g_{i}}(\delta_{i}s)|\le O(\delta_{i}^{4}|s|^{4})$,
we get
\begin{equation}
|I_{4}(u_{i},r)|=\frac{n-2}{2}\delta_{i}\int\limits _{\partial'B_{R_{i}}^{+}}\left(\bar{s}^{k}\partial_{k}v_{i}+\frac{n-2}{2}v_{i}\right)h_{g_{i}}(\delta_{i}s)v_{i}d\bar{s}\le c\delta_{i}^{5-2\eta}.\label{eq:poho3}
\end{equation}
Using the expansion of the metric, it easy to check that
\begin{equation}
|I_{3}(u_{i},r)|\le c\delta_{i}^{2-2\eta}.\label{eq:poho2}
\end{equation}
Finally, by Claim 1 of Proposition \ref{prop:4.1}, by (\ref{eq:doppioint})
and by (\ref{eq:Iam}) we get
\begin{align*}
\lim_{i\rightarrow\infty} & \int\limits _{\partial'B_{R_{i}}^{+}}\left(\bar{s}^{k}\partial_{k}v_{i}+\frac{n-2}{2}v_{i}\right)\beta_{i}(\delta_{i}s)v_{i}d\bar{s}\\
 & =\beta(x_{0})\int\limits _{\mathbb{R}^{n-1}}\left(\bar{s}^{k}\partial_{k}U+\frac{n-2}{2}U\right)Ud\bar{s}\\
 & =\frac{n-2}{2}\beta(x_{0})\int\limits _{\mathbb{R}^{n-1}}\frac{1-|\bar{s}|^{2}}{\left[1+|\bar{s}|^{2}\right]^{n-1}}d\bar{s}\\
 & =\frac{n-2}{2}\omega_{n-2}\beta(x_{0})\left[I_{n-1}^{n-2}-I_{n-1}^{n}\right]=-\frac{n-2}{n-1}\beta(x_{0})\omega_{n-2}I_{n-1}^{n}=:B>0,
\end{align*}
so 
\begin{equation}
I_{6}(u_{i},r)=\varepsilon_{2,i}\delta_{i}(B+o(1)).\label{eq:poho4}
\end{equation}
In a similar way we proceed for $I_{5}$. In fact, by (\ref{eq:doppioint}),
(\ref{eq:Iam}) and (\ref{eq:t-integrali}), we have
\begin{multline}
\lim_{i\rightarrow\infty}\int\limits _{B_{R_{i}}^{+}}\left(\bar{s}^{a}\partial_{a}v_{i}+\frac{n-2}{2}v_{i}\right)\alpha_{i}(\delta_{i}s)v_{i}d\bar{s}=\\
=\alpha(x_{0})\int\limits _{\mathbb{R}_{+}^{n}}\left(\bar{s}^{a}\partial_{a}U+\frac{n-2}{2}U\right)Ud\bar{s}\\
=\frac{n-2}{2}\alpha(x_{0})\int\limits _{\mathbb{R}^{+}}\int\limits _{\mathbb{R}^{n-1}}\frac{1-|s_{n}|^{2}-|\bar{s}|^{2}}{\left[|1+s_{n}|^{2}+|\bar{s}|^{2}\right]^{n-1}}d\bar{s}ds_{n}\\
=\frac{n-2}{2}\alpha(x_{0})\omega_{n-2}\left[I_{n-1}^{n-2}\int_{0}^{\infty}\frac{1-t^{2}}{(1+t)^{n-1}}dt-I_{n-1}^{n}\int_{0}^{\infty}\frac{1}{(1+t)^{n-3}}dt\right]\\
=\frac{n-2}{2}\alpha(x_{0})\omega_{n-2}\left[I_{n-1}^{n-2}\frac{n-5}{(n-3)(n-4)}-I_{n-1}^{n}\frac{1}{n-4}\right]\\
=-\frac{2(n-2)}{(n-1)(n-4)}\alpha(x_{0})I_{n-1}^{n}\omega_{n-2}=:A>0\label{eq:pohosign}
\end{multline}
and thus 
\begin{equation}
I_{5}(u_{i},r)=\varepsilon_{1,i}\delta_{i}^{2}(A+o(1)).\label{eq:poho5}
\end{equation}

Concluding, by (\ref{eq:poho1}), (\ref{eq:poho2}), (\ref{eq:poho3}),
(\ref{eq:poho4}), (\ref{eq:poho4}). we get 
\[
-c\delta_{i}^{2-2\eta}+(A+o(1))\varepsilon_{1,i}\delta_{i}^{2}+(B+o(1))\varepsilon_{2,i}\delta_{i}\le\delta_{i}^{\lambda_{i}(n-2)}
\]
 which is possible only if $\varepsilon_{1,i}\delta_{i},\varepsilon_{2,i}\rightarrow0$.
\end{proof}
Since $\varepsilon_{1,i}\delta_{i},\varepsilon_{2,i}\rightarrow0$
by Prop. \ref{prop:eps-i}, the proof of the next proposition is analogous
to Prop. 4.3 in \cite{Al}
\begin{prop}
\label{prop:4.3}Let $x_{i}\rightarrow x_{0}$ be an isolated simple
blow-up point for $\left\{ u_{i}\right\} _{i}$ and $\alpha,\beta<0$.
Then there exist $C,\rho>0$ such that 
\begin{enumerate}
\item $M_{i}u_{i}(\psi_{i}(y))\le C|y|^{2-n}$ for all $y\in B_{\rho}^{+}(0)\smallsetminus\left\{ 0\right\} $;
\item $M_{i}u_{i}(\psi_{i}(y))\ge C^{-1}G_{i}(y)$ for all $y\in B_{\rho}^{+}(0)\smallsetminus B_{r_{i}}^{+}(0)$
where $r_{i}:=R_{i}M_{i}^{\frac{2}{2-n}}$ and $G_{i}$ is the Green\textquoteright s
function which solves
\[
\left\{ \begin{array}{ccc}
L_{g_{i}}G_{i}=0 &  & \text{in }B_{\rho}^{+}(0)\smallsetminus\left\{ 0\right\} \\
G_{i}=0 &  & \text{on }\partial^{+}B_{\rho}^{+}(0)\\
B_{g_{i}}G_{i}=0 &  & \text{on }\partial'B_{\rho}^{+}(0)\smallsetminus\left\{ 0\right\} 
\end{array}\right.
\]
\end{enumerate}
and $|y|^{n-2}G_{i}(y)\rightarrow1$ as $|z|\rightarrow0$.
\end{prop}

By Proposition \ref{prop:4.1} and Proposition \ref{prop:4.3} we
have that, if $x_{i}\rightarrow x_{0}$ is an isolated simple blow-up
point for $\left\{ u_{i}\right\} _{i}$, then it holds
\[
v_{i}\le CU\text{ in }B_{\rho M_{i}^{\frac{2}{2-n}}}^{+}(0).
\]

\subsection{Blowup estimates\label{sec:Blowup-estimates}}

Our aim is to provide a fine estimate for the approximation of the
rescaled solution near an isolated simple blow up point. In this section
$x_{i}\rightarrow x_{0}$ is an isolated simple blowup point for a
sequence $\left\{ u_{i}\right\} _{i}$ of solutions of (\ref{eq:Prob-i}).
We will work in the conformal Fermi coordinates in a neighborhood
of $x_{i}$.

Set $\tilde{u}_{i}=\Lambda_{x_{i}}^{-1}u_{i}$ and
\begin{equation}
\delta_{i}:=\tilde{u}_{i}^{\frac{2}{2-n}}(x_{i})=u_{i}^{\frac{2}{2-n}}(x_{i})=M_{i}^{\frac{2}{2-n}}\ \ \ v_{i}(y):=\delta_{i}^{\frac{n-2}{2}}u_{i}(\delta_{i}y)\text{ for }y\in B_{\frac{R}{\delta_{i}}}^{+}(0).\label{eq:deltai}
\end{equation}
 Then $v_{i}$ satisfies 
\begin{equation}
\left\{ \begin{array}{cc}
L_{\hat{g}_{i}}v_{i}-\varepsilon_{1,i}\alpha_{i}(\delta_{i}y)v_{i}=0 & \text{ in }B_{\frac{R}{\delta_{i}}}^{+}(0)\\
B_{\hat{g}_{i}}v_{i}+(n-2)v_{i}^{\frac{n}{n-2}}-\varepsilon_{2,i}\beta_{i}(\delta_{i}y)v_{i}=0 & \text{ on }\partial'B_{\frac{R}{\delta_{i}}}^{+}(0)
\end{array}\right.\label{eq:Prob-hat}
\end{equation}
 where $\hat{g}_{i}:=\tilde{g}_{i}(\delta_{i}y)=\Lambda_{x_{i}}^{\frac{4}{n-2}}(\delta_{i}y)g(\delta_{i}y)$,
and $\alpha_{i}(y)=\Lambda_{x_{i}}^{-\frac{2}{n-2}}(y)\alpha(y)$,
$\beta_{i}(y)=\Lambda_{x_{i}}^{-\frac{2}{n-2}}(y)\beta(y)$. 

The estimates that follow are similar to the ones of \cite[Lemma 6.1]{Al},
\cite[Section 4]{GM20} and \cite[Section 5]{GMdcds}, where the main
differences concern the terms which contain the linear perturbations.
\begin{lem}
\label{lem:coreLemma}Assume $n\ge8$. Let $\gamma_{x_{i}}$ be defined
in (\ref{eq:vqdef}). There exist $R,C>0$ such that 
\[
|v_{i}(y)-U(y)-\delta_{i}^{2}\gamma_{x_{i}}(y)|\le C\left(\delta_{i}^{3}+\varepsilon_{1,i}\delta_{i}^{2}+\varepsilon_{2,i}\delta_{i}\right)
\]
for $|y|\le R/\delta_{i}$.
\end{lem}

\begin{proof}
Let $y_{i}$ such that 
\[
\mu_{i}:=\max_{|y|\le R/\delta_{i}}|v_{i}(y)-U(y)-\delta_{i}^{2}\gamma_{x_{i}}(y)|=|v_{i}(y_{i})-U(y_{i})-\delta_{i}^{2}\gamma_{x_{i}}(y_{i})|.
\]
We can assume, without loss of generality, that $|y_{i}|\le\frac{R}{2\delta_{i}}.$

In fact, suppose that there exists $c>0$ such that $|y_{i}|>\frac{c}{\delta_{i}}$
for all $i$. Then, since $v_{i}(y)\le CU(y)$, and by (\ref{eq:gradvq}),
we get the inequality
\[
|v_{i}(y_{i})-U(y_{i})-\delta_{i}^{2}\gamma_{x_{i}}(y_{i})|\le C\left(|y_{i}|^{2-n}+\delta_{i}^{2}|y_{i}|^{4-n}\right)\le C\delta_{i}^{n-2}
\]
which proves the Lemma. So, in the next we will suppose $|y_{i}|\le\frac{R}{2\delta_{i}}$.
This will be useful later.

By contradiction, suppose that 
\begin{equation}
\max\left\{ \mu_{i}^{-1}\delta_{i}^{3},\mu_{i}^{-1}\varepsilon_{1,i}\delta_{i}^{2},\mu_{i}^{-1}\varepsilon_{2,i}\delta_{i}\right\} \rightarrow0\text{ when }i\rightarrow\infty.\label{eq:ipass}
\end{equation}
Defined 
\[
w_{i}(y):=\mu_{i}^{-1}\left(v_{i}(y)-U(y)-\delta_{i}^{2}\gamma_{x_{i}}(y)\right)\text{ for }|y|\le R/\delta_{i},
\]
we have, by direct computation, that $w_{i}$ satisfies 
\begin{equation}
\left\{ \begin{array}{cc}
L_{\hat{g}_{i}}w_{i}=A_{i} & \text{ in }B_{\frac{R}{\delta_{i}}}^{+}(0)\\
B_{\hat{g}_{i}}w_{i}+b_{i}w_{i}=F_{i} & \text{ on }\partial'B_{\frac{R}{\delta_{i}}}^{+}(0)
\end{array}\right.\label{eq:wi}
\end{equation}
where 
\begin{align*}
b_{i}= & (n-2)\frac{v_{i}^{\frac{n}{n-2}}-(U+\delta_{i}^{2}\gamma_{x_{i}})^{\frac{n}{n-2}}}{v_{i}-U-\delta_{i}^{2}\gamma_{x_{i}}},\\
Q_{i}= & -\frac{1}{\mu_{i}}\left\{ \left(L_{\hat{g}_{i}}-\Delta\right)(U+\delta_{i}^{2}\gamma_{x_{i}})+\delta_{i}^{2}\Delta\gamma_{x_{i}}\right\} ,\\
A_{i}= & Q_{i}+\frac{\varepsilon_{1,i}\delta_{i}^{2}}{\mu_{i}}\alpha_{i}(\delta_{i}y)v_{i}(y),\\
\bar{Q}_{i}= & -\frac{1}{\mu_{i}}\left\{ (n-2)(U+\delta_{i}^{2}\gamma_{x_{i}})^{\frac{n}{n-2}}-(n-2)U^{\frac{n}{n-2}}-n\delta_{i}^{2}U^{\frac{2}{n-2}}\gamma_{x_{i}}-\frac{n-2}{2}h_{\hat{g}_{i}}(U+\delta_{i}^{2}\gamma_{x_{i}})\right\} ,\\
F_{i}= & \bar{Q}_{i}+\frac{\varepsilon_{2,i}\delta_{i}}{\mu_{i}}\beta_{i}(\delta_{i}y)v_{i}(y).
\end{align*}
We will estimate terms $b_{i},A_{i,}F_{i}$ obtaining that the sequence
$w_{i}$ converges in $C_{\text{loc}}^{2}(\mathbb{R}_{+}^{n})$ to
some $w$ solution of 
\begin{equation}
\left\{ \begin{array}{cc}
\Delta w=0 & \text{ in }\mathbb{R}_{+}^{n}\\
\frac{\partial}{\partial\nu}w+nU^{\frac{n}{n-2}}w=0 & \text{ on }\partial\mathbb{R}_{+}^{n}
\end{array}\right.,\label{eq:diff-w}
\end{equation}
then we will derive a contradiction using (\ref{eq:ipass}). 

Since $v_{i}\rightarrow U$ in $C_{\text{loc}}^{2}(\mathbb{R}_{+}^{n})$
we have, at once, 
\begin{align}
b_{i} & \rightarrow nU^{\frac{2}{n-2}}\text{ in }C_{\text{loc}}^{2}(\mathbb{R}_{+}^{n});\label{eq:b1}\\
|b_{i}(y)| & \le(1+|y|)^{-2}\text{ for }|y|\le R/\delta_{i}.\label{eq:b2}
\end{align}
We proceed now by estimating $Q_{i}$ and $\bar{Q}_{i}$. As in \cite[Lemma 11]{GMdcds},
using the expansion of the metric and the decays properties of $U$
and $\gamma_{x_{i}}$ we obtain
\begin{equation}
Q_{i}=O(\mu_{i}^{-1}\delta_{i}^{3}\left(1+|y|\right)^{3-n})\label{eq:Q}
\end{equation}
and 
\begin{equation}
\bar{Q}_{i}=O(\mu_{i}^{-1}\delta_{i}^{4}\left(1+|y|\right)^{5-n}).\label{eq:Qbar}
\end{equation}
Since $|v_{i}(y)|\le CU(y)$ from (\ref{eq:Q}) and (\ref{eq:Qbar})
we get

\begin{align}
A_{i} & =O(\mu_{i}^{-1}\delta_{i}^{3}\left(1+|y|\right)^{3-n})+O(\mu_{i}^{-1}\varepsilon_{1,i}\delta_{i}^{2}\left(1+|y|\right)^{2-n}),\label{eq:A-B}\\
F_{i} & =O(\mu_{i}^{-1}\delta_{i}^{4}\left(1+|y|\right)^{5-n})+O(\mu_{i}^{-1}\varepsilon_{2,i}\delta_{i}\left(1+|y|\right)^{2-n}).\nonumber 
\end{align}
In light of (\ref{eq:ipass}) we also have $A_{i}\in L^{p}(B_{R/\delta_{i}}^{+})$
and $F_{i}\in L^{p}(\partial'B_{R/\delta_{i}}^{+})$ for all $p\ge2$. 

Finally we remark that $|w_{i}(y)|\le1$, so by (\ref{eq:ipass})
(\ref{eq:b1}), (\ref{eq:b2}), (\ref{eq:A-B}) and by standard elliptic
estimates we conclude that, up to subsequence, $\left\{ w_{i}\right\} _{i}$
converges in $C_{\text{loc}}^{2}(\mathbb{R}_{+}^{n})$ to some $w$
solution of (\ref{eq:diff-w}) as claimed at the beginning of the
proof.

The next step is to prove that $|w(y)|\le C(1+|y|^{-1})$ for $y\in\mathbb{R}_{+}^{n}$.
Consider $G_{i}$ the Green function for the conformal Laplacian $L_{\hat{g}_{i}}$
defined on $B_{r/\delta_{i}}^{+}$ with boundary conditions $B_{\hat{g}_{i}}G_{i}=0$
on $\partial'B_{r/\delta_{i}}^{+}$ and $G_{i}=0$ on $\partial^{+}B_{r/\delta_{i}}^{+}$.
It is well known that $G_{i}=O(|\xi-y|^{2-n})$. By the Green formula
and by (\ref{eq:A-B}) we have
\begin{align*}
w_{i}(y)= & -\int_{B_{\frac{R}{\delta_{i}}}^{+}}G_{i}(\xi,y)A_{i}(\xi)d\mu_{\hat{g}_{i}}(\xi)-\int_{\partial^{+}B_{\frac{R}{\delta_{i}}}^{+}}\frac{\partial G_{i}}{\partial\nu}(\xi,y)w_{i}(\xi)d\sigma_{\hat{g}_{i}}(\xi)\\
 & +\int_{\partial'B_{\frac{R}{\delta_{i}}}^{+}}G_{i}(\xi,y)\left(b_{i}(\xi)w_{i}(\xi)-F_{i}(\xi)\right)d\sigma_{\hat{g}_{i}}(\xi),
\end{align*}
so 
\begin{align*}
|w_{i}(y)| & \le\frac{\delta_{i}^{3}}{\mu_{i}}\int_{B_{\frac{R}{\delta_{i}}}^{+}}|\xi-y|^{2-n}(1+|\xi|)^{3-n}d\xi+\frac{\varepsilon_{1,i}\delta_{i}^{2}}{\mu_{i}}\int_{B_{\frac{R}{\delta_{i}}}^{+}}|\xi-y|^{2-n}(1+|\xi|)^{2-n}d\xi\\
 & +\int_{\partial^{+}B_{\frac{R}{\delta_{i}}}^{+}}|\xi-y|^{1-n}w_{i}(\xi)d\sigma(\xi)\\
 & +\int_{\partial'B_{\frac{R}{\delta_{i}}}^{+}}|\bar{\xi}-y|^{2-n}\left((1+|\bar{\xi}|)^{-2}+\frac{\delta_{i}^{4}}{\mu_{i}}(1+|\bar{\xi}|)^{5-n}+\frac{\varepsilon_{2,i}\delta_{i}}{\mu_{i}}(1+|\bar{\xi}|)^{2-n}\right)d\bar{\xi}.
\end{align*}
Notice that in the third integral we used that $|y|\le\frac{R}{2\delta_{i}}$
to estimate $|\xi-y|\ge|\xi|-|y|\ge\frac{R}{2\delta_{i}}$ on $\partial^{+}B_{R/\delta_{i}}^{+}$.
Moreover, since $v_{i}(\xi)\le CU(\xi)$, we get
\begin{equation}
|w_{i}(\xi)|\le\frac{C}{\mu_{i}}\left(\left(1+|\xi|\right)^{2-n}+\delta_{i}^{2}\left(1+|\xi|\right)^{4-n}\right)\le C\frac{\delta_{i}^{n-2}}{\mu_{i}}\text{ on }\partial^{+}B_{R/\delta_{i}}^{+};\label{eq:wibordo}
\end{equation}
hence 
\begin{equation}
\int_{\partial^{+}B_{\frac{R}{\delta_{i}}}^{+}}|\xi-y|^{1-n}w_{i}(\xi)d\sigma(\xi)\le C\int_{\partial^{+}B_{\frac{R}{\delta_{i}}}^{+}}\frac{\delta_{i}^{2n-3}}{\mu_{i}}d\sigma_{\hat{g}_{i}}(\xi)\le C\frac{\delta_{i}^{n-2}}{\mu_{i}}.\label{eq:stimaW1}
\end{equation}
 For the other terms we use the formula 
\begin{equation}
\int_{\mathbb{R}^{m}}|\xi-y|^{l-m}(1+|\xi|)^{-\eta}d\xi\le C(1+|y|)^{l-\eta}\label{eq:ALstimagreen}
\end{equation}
where $y\in\mathbb{R}^{m+k}\supseteq\mathbb{R}^{m}$, $\eta,l\in\mathbb{N}$,
$0<l<\eta<m$ (see \cite[Lemma 9.2]{Al} and \cite{Au,Gi}) . We get
\begin{equation}
\frac{\delta_{i}^{3}}{\mu_{i}}\int_{B_{\frac{R}{\delta_{i}}}^{+}}|\xi-y|^{2-n}(1+|\xi|)^{3-n}d\xi\le C\frac{\delta_{i}^{3}}{\mu_{i}}(1+|y|)^{5-n},\label{eq:stimaW2}
\end{equation}
\begin{equation}
\frac{\varepsilon_{1,i}\delta_{i}^{2}}{\mu_{i}}\int_{B_{\frac{R}{\delta_{i}}}^{+}}|\xi-y|^{2-n}(1+|\xi|)^{2-n}d\xi\le C\frac{\varepsilon_{1,i}\delta_{i}^{2}}{\mu_{i}}(1+|y|)^{4-n},\label{eq:stimaW2-bis}
\end{equation}

\begin{equation}
\int_{\partial'B_{\frac{R}{\delta_{i}}}^{+}}|\bar{\xi}-y|^{2-n}(1+|\bar{\xi}|)^{-2}d\bar{\xi}\le(1+|y|)^{-1},\label{eq:stimaW3}
\end{equation}
\begin{equation}
\frac{\delta_{i}^{4}}{\mu_{i}}\int_{\partial'B_{\frac{R}{\delta_{i}}}^{+}}|\bar{\xi}-y|^{2-n}(1+|\bar{\xi}|)^{5-n}d\bar{\xi}\le C\frac{\delta_{i}^{4}}{\mu_{i}}(1+|y|)^{6-n},\label{eq:stimaW4}
\end{equation}

\begin{equation}
\frac{\varepsilon_{2,i}\delta_{i}}{\mu_{i}}\int_{\partial'B_{\frac{R}{\delta_{i}}}^{+}}|\bar{\xi}-y|^{2-n}(1+|\bar{\xi}|)^{2-n}d\bar{\xi}\le C\frac{\varepsilon_{2,i}\delta_{i}}{\mu_{i}}(1+|y|)^{3-n}.\label{eq:stimaW5}
\end{equation}
By the previous estimates we infer that, for $|y|\le\frac{R}{2\delta_{i}}$,
\[
|w_{i}(y)|\le C\left((1+|y|)^{-1}+\frac{\delta_{i}^{3}}{\mu_{i}}(1+|y|)^{5-n}+\frac{\varepsilon_{1,i}\delta_{i}^{2}}{\mu_{i}}(1+|y|)^{4-n}+\frac{\varepsilon_{2,i}\delta_{i}}{\mu_{i}}(1+|y|)^{3-n}\right)
\]
so by assumption (\ref{eq:ipass}) we prove 
\begin{equation}
|w(y)|\le C(1+|y|)^{-1}\text{ for }y\in\mathbb{R}_{+}^{n}\label{eq:stimaWass}
\end{equation}
as claimed.

Finally we notice that, since $v_{i}\rightarrow U$ near $0$, and,
by (\ref{eq:dervq}), we have $w_{i}(0)\rightarrow0$ as well as $\frac{\partial w_{i}}{\partial y_{j}}(0)\rightarrow0$
for $j=1,\dots,n-1$. This implies that 
\begin{equation}
w(0)=\frac{\partial w}{\partial y_{1}}(0)=\dots=\frac{\partial w}{\partial y_{n-1}}(0)=0.\label{eq:W(0)}
\end{equation}
We are ready now to prove the contradiction. In fact, it is known
(see \cite[Lemma 2]{Al}) that any solution of (\ref{eq:diff-w})
that decays as (\ref{eq:stimaWass}) is a linear combination of $\frac{\partial U}{\partial y_{1}},\dots,\frac{\partial U}{\partial y_{n-1}},\frac{n-2}{2}U+y^{b}\frac{\partial U}{\partial y_{b}}$.
This fact, combined with (\ref{eq:W(0)}), implies that $w\equiv0$. 

Now, on the one hand $|y_{i}|\le\frac{R}{2\delta_{i}}$, so estimate
(\ref{eq:stimaWass}) holds; on the other hand, since $w_{i}(y_{i})=1$
and $w\equiv0$, we get $|y_{i}|\rightarrow\infty$, obtaining
\[
1=w_{i}(y_{i})\le C(1+|y_{i}|)^{-1}\rightarrow0
\]
 which gives us the contradiction. 
\end{proof}
\begin{lem}
\label{lem:taui}Assume $n\ge8$ and $\alpha,\beta<0$. There exists
$R,C>0$ such that 
\[
\varepsilon_{1,i}\delta_{i}^{2}+\varepsilon_{2,i}\delta_{i}\le C\delta_{i}^{3}
\]
for $|y|\le R/\delta_{i}$.
\end{lem}

\begin{proof}
We proceed by contradiction, supposing that 
\begin{equation}
\left(\varepsilon_{1,i}\delta_{i}^{2}+\varepsilon_{2,i}\delta_{i}\right)^{-1}\delta_{i}^{3}\rightarrow0\text{ when }i\rightarrow\infty.\label{eq:ipasstau}
\end{equation}
Thus, by Lemma \ref{lem:coreLemma}, we have 
\[
|v_{i}(y)-U(y)-\delta_{i}^{2}\gamma_{x_{i}}(y)|\le C(\varepsilon_{1,i}\delta_{i}^{2}+\varepsilon_{2,i}\delta_{i})\text{ for }|y|\le R/\delta_{i}.
\]
We define, similarly to Lemma \ref{lem:coreLemma},
\[
w_{i}(y):=\frac{1}{\varepsilon_{1,i}\delta_{i}^{2}+\varepsilon_{2,i}\delta_{i}}\left(v_{i}(y)-U(y)-\delta_{i}^{2}\gamma_{x_{i}}(y)\right)\text{ for }|y|\le R/\delta_{i},
\]
and we have that $w_{i}$ satisfies (\ref{eq:wi}) where $b_{i}$
is as in Lemma \ref{lem:coreLemma} and
\begin{align*}
Q_{i}= & -\frac{1}{\varepsilon_{1,i}\delta_{i}^{2}+\varepsilon_{2,i}\delta_{i}}\left\{ \left(L_{\hat{g}_{i}}-\Delta\right)(U+\delta_{i}^{2}\gamma_{x_{i}})+\delta_{i}^{2}\Delta\gamma_{x_{i}}\right\} ,\\
A_{i}= & Q_{i}+\frac{\varepsilon_{1,i}\delta_{i}^{2}}{\varepsilon_{1,i}\delta_{i}^{2}+\varepsilon_{2,i}\delta_{i}}\alpha_{i}(\delta_{i}y)v_{i}(y),\\
\bar{Q}_{i}= & -\frac{1}{\varepsilon_{1,i}\delta_{i}^{2}+\varepsilon_{2,i}\delta_{i}}\left\{ (n-2)(U+\delta_{i}^{2}\gamma_{x_{i}})^{\frac{n}{n-2}}-(n-2)U^{\frac{n}{n-2}}-n\delta_{i}^{2}U^{\frac{2}{n-2}}\gamma_{x_{i}}-\frac{n-2}{2}h_{\hat{g}_{i}}(U+\delta_{i}^{2}\gamma_{x_{i}})\right\} ,\\
F_{i}= & \bar{Q}_{i}+\frac{\varepsilon_{2,i}\delta_{i}}{\varepsilon_{1,i}\delta_{i}^{2}+\varepsilon_{2,i}\delta_{i}}\beta_{i}(\delta_{i}y)v_{i}(y).
\end{align*}
As before, $b_{i}$ satisfies inequality (\ref{eq:b2}) while 

\begin{align}
A_{i} & =O\left(\frac{\delta_{i}^{3}}{\varepsilon_{1,i}\delta_{i}^{2}+\varepsilon_{2,i}\delta_{i}}\left(1+|y|\right)^{3-n}\right)+O\left(\frac{\varepsilon_{1,i}\delta_{i}^{2}}{\varepsilon_{1,i}\delta_{i}^{2}+\varepsilon_{2,i}\delta_{i}}\left(1+|y|\right)^{2-n}\right)\label{eq:A-B-1}\\
F_{i} & =O\left(\frac{\delta_{i}^{4}}{\varepsilon_{1,i}\delta_{i}^{2}+\varepsilon_{2,i}\delta_{i}}\left(1+|y|\right)^{5-n}\right)+O\left(\frac{\varepsilon_{2,i}\delta_{i}}{\varepsilon_{1,i}\delta_{i}^{2}+\varepsilon_{2,i}\delta_{i}}\left(1+|y|\right)^{2-n}\right),\nonumber 
\end{align}
so by classic elliptic estimates we can prove that the sequence $w_{i}$
converges in $C_{\text{loc}}^{2}(\mathbb{R}_{+}^{n})$ to some $w$.

We proceed as in Lemma \ref{lem:coreLemma} to deduce that, by (\ref{eq:ipasstau})
and since $\frac{\varepsilon_{1,i}\delta_{i}^{2}}{\varepsilon_{1,i}\delta_{i}^{2}+\varepsilon_{2,i}\delta_{i}}\le1$,
$\frac{\varepsilon_{2,i}\delta_{i}}{\varepsilon_{1,i}\delta_{i}^{2}+\varepsilon_{2,i}\delta_{i}}\le1$,
\begin{align}
|w_{i}(y)| & \le C\left((1+|y|)^{-1}+\frac{\delta_{i}^{3}(1+|y|)^{5-n}}{\varepsilon_{1,i}\delta_{i}^{2}+\varepsilon_{2,i}\delta_{i}}+\frac{\varepsilon_{1,i}\delta_{i}^{2}\left(1+|y|\right)^{4-n}}{\varepsilon_{1,i}\delta_{i}^{2}+\varepsilon_{2,i}\delta_{i}}+\frac{\varepsilon_{2,i}\delta_{i}\left(1+|y|\right)^{3-n}}{\varepsilon_{1,i}\delta_{i}^{2}+\varepsilon_{2,i}\delta_{i}}\right)\nonumber \\
 & \le C\left((1+|y|)^{-1}\right)\text{ for }|y|\le\frac{R}{2\delta_{i}}.\label{eq:decWtau}
\end{align}

Now let $j_{n}$ defined as in (\ref{eq:jn}). Indeed, since $w_{i}$
satisfies (\ref{eq:wi}), integrating by parts we obtain
\begin{multline}
\int_{\partial'B_{\frac{R}{\delta_{i}}}^{+}}j_{n}F_{i}d\sigma_{\hat{g}_{i}}=\int_{\partial'B_{\frac{R}{\delta_{i}}}^{+}}j_{n}\left[B_{\hat{g}_{i}}w_{i}+b_{i}w_{i}\right]d\sigma_{\hat{g}_{i}}\\
=\int_{\partial'B_{\frac{R}{\delta_{i}}}^{+}}w_{i}\left[B_{\hat{g}_{i}}j_{n}+b_{i}j_{n}\right]d\sigma_{\hat{g}_{i}}+\int_{\partial^{+}B_{\frac{R}{\delta_{i}}}^{+}}\left[\frac{\partial j_{n}}{\partial\eta_{i}}w_{i}-\frac{\partial w_{i}}{\partial\eta_{i}}j_{n}\right]d\sigma_{\hat{g}_{i}}\\
+\int_{B_{\frac{R}{\delta_{i}}}^{+}}\left[w_{i}L_{\hat{g}_{i}}j_{n}-j_{n}L_{\hat{g}_{i}}w_{i}\right]d\mu_{\hat{g}_{i}}\label{eq:parts}
\end{multline}
where $\eta_{i}$ is the inward unit normal vector to $\partial^{+}B_{\frac{R}{\delta_{i}}}^{+}$.
One can check easily that

\[
\lim_{i\rightarrow+\infty}\int_{\partial'B_{\frac{R}{\delta_{i}}}^{+}}j_{n}\bar{Q}_{i}d\sigma_{\hat{g}_{i}}=0.
\]
Also, since $\beta_{i}(\delta_{i}y)=\Lambda_{x_{i}}^{-\frac{2}{n-2}}(\delta_{i}y)\beta(\delta_{i}y)$,
and $\beta<0$, by Proposition \ref{prop:4.1}, we have 
\[
\beta_{i}(\delta_{i}y)v_{i}(y)\rightarrow\beta(x_{0})U(y)\text{ for }i\rightarrow+\infty.
\]
and thus, as in Proposition \ref{prop:eps-i}
\begin{equation}
\lim_{i\rightarrow+\infty}\int_{\partial'B_{\frac{R}{\delta_{i}}}^{+}}\beta_{i}(\delta_{i}y)v_{i}(y)j_{n}(y)=\frac{n-2}{2}\beta(x_{0})\int_{\mathbb{R}^{n-1}}\frac{1-|\bar{y}|^{2}}{\left(1+|\bar{y}|^{2}\right)^{n-1}}=:B>0\label{eq:eq:avj}
\end{equation}
so
\begin{equation}
\int_{\partial'B_{\frac{R}{\delta_{i}}}^{+}}j_{n}F_{i}d\sigma_{\hat{g}_{i}}=\frac{\varepsilon_{2,i}\delta_{i}}{\varepsilon_{1,i}\delta_{i}^{2}+\varepsilon_{2,i}\delta_{i}}(B+o(1)).\label{eq:neg}
\end{equation}
By (\ref{eq:parts}) and (\ref{eq:neg}) we derive a contradiction.
Indeed, by the decay of $j_{n}$ and by the decay of $w_{i}$, given
by (\ref{eq:decWtau}) and by (\ref{eq:ipasstau}), we have
\begin{equation}
\lim_{i\rightarrow+\infty}\int_{\partial^{+}B_{\frac{R}{\delta_{i}}}^{+}}\left[\frac{\partial j_{n}}{\partial\eta_{i}}w_{i}-\frac{\partial w_{i}}{\partial\eta_{i}}j_{n}\right]d\sigma_{\hat{g}_{i}}=0.\label{eq:nullo1}
\end{equation}
Since $\Delta j_{n}=0$, one can check that 
\begin{equation}
\lim_{i\rightarrow+\infty}\int_{B_{\frac{R}{\delta_{i}}}^{+}}w_{i}L_{\hat{g}_{i}}j_{n}d\mu_{\hat{g}_{i}}=0.\label{eq:nullo3}
\end{equation}
Also, we can prove that
\begin{equation}
\lim_{i\rightarrow+\infty}\int_{B_{\frac{R}{\delta_{i}}}^{+}}j_{n}Q_{i}d\mu_{\hat{g}_{i}}=0.\label{eq:nullo2}
\end{equation}
Finally 
\begin{align}
\lim_{i\rightarrow+\infty}\int_{\partial'B_{\frac{R}{\delta_{i}}}^{+}}w_{i}\left[B_{\hat{g}_{i}}j_{n}+b_{i}j_{n}\right]d\sigma_{\hat{g}_{i}} & =\int_{\partial\mathbb{R}_{+}^{n}}w\left[\frac{\partial j_{n}}{\partial y_{n}}+nU^{\frac{2}{n-2}}j_{n}\right]d\sigma_{\hat{g}_{i}}=0\label{eq:nullofinale}
\end{align}
since $\frac{\partial j_{n}}{\partial y_{n}}+nU^{\frac{2}{n-2}}j_{n}=0$
when $y_{n}=0$. 

In light of (\ref{eq:nullo1}) (\ref{eq:nullo2}) and (\ref{eq:nullo3})
we infer, by (\ref{eq:parts}), that 
\begin{align}
\int_{\partial'B_{\frac{R}{\delta_{i}}}^{+}}j_{n}F_{i}d\sigma_{\hat{g}_{i}} & =-\int_{B_{\frac{R}{\delta_{i}}}^{+}}\left[j_{n}A_{i}w_{i}\right]d\mu_{\hat{g}_{i}}+o(1).\label{eq:parts-1}\\
 & =-\frac{\varepsilon_{1,i}\delta_{i}^{2}}{\varepsilon_{1,i}\delta_{i}^{2}+\varepsilon_{2,i}\delta_{i}}\int_{B_{\frac{R}{\delta_{i}}}^{+}}j_{n}(y)\alpha_{i}(\delta_{i}y)v_{i}(y)d\mu_{\hat{g}_{i}}+o(1).\nonumber 
\end{align}
Again we have $\alpha_{i}(\delta_{i}y)v_{i}(y)\rightarrow\alpha(x_{0})U(y)\text{ for }i\rightarrow+\infty$
and $\alpha<0$, so, proceeding as in Proposition \ref{prop:eps-i},
we have 
\begin{equation}
\lim_{i\rightarrow\infty}\int_{B_{\frac{R}{\delta_{i}}}^{+}}j_{n}(y)\alpha_{i}(\delta_{i}y)v_{i}(y)d\mu_{\hat{g}_{i}}=\alpha(x_{0})\lim_{i\rightarrow\infty}\int\limits _{\mathbb{R}_{+}^{n}}\left(s^{a}\partial_{a}v_{i}+\frac{n-2}{2}v_{i}\right)v_{i}ds=:A>0,\label{eq:bvj}
\end{equation}
so
\begin{equation}
\int_{B_{\frac{R}{\delta_{i}}}^{+}}\left[j_{n}A_{i}w_{i}\right]=-\frac{\varepsilon_{1,i}\delta_{i}^{2}}{\varepsilon_{1,i}\delta_{i}^{2}+\varepsilon_{2,i}\delta_{i}}(A+o(1))\label{eq:bneg}
\end{equation}
and, by (\ref{eq:parts-1}), (\ref{eq:neg}) and (\ref{eq:bneg}),
we have 
\begin{equation}
\frac{\varepsilon_{2,i}\delta_{i}}{\varepsilon_{1,i}\delta_{i}^{2}+\varepsilon_{2,i}\delta_{i}}(B+o(1))=-\frac{\varepsilon_{1,i}\delta_{i}^{2}}{\varepsilon_{1,i}\delta_{i}^{2}+\varepsilon_{2,i}\delta_{i}}(A+o(1)).\label{eq:contra}
\end{equation}
Since $\frac{\varepsilon_{2,i}\delta_{i}}{\varepsilon_{1,i}\delta_{i}^{2}+\varepsilon_{2,i}\delta_{i}}$
and $\frac{\varepsilon_{1,i}\delta_{i}^{2}}{\varepsilon_{1,i}\delta_{i}^{2}+\varepsilon_{2,i}\delta_{i}}$
cannot vanish simultaneously while $i\rightarrow\infty$, equation
(\ref{eq:contra}) leads us to a contradiction. 
\end{proof}
The above lemmas are the core of the following proposition, in which
we iterate the procedure of Lemma \ref{lem:coreLemma}, to obtain
better estimates of the rescaled solution $v_{i}$ of (\ref{eq:Prob-hat})
around the isolated simple blow up point $x_{i}\rightarrow x_{0}$.
\begin{prop}
\label{prop:stimawi}Assume $n\ge8$ and $\alpha,\beta<0$. Let $\gamma_{x_{i}}$
be defined in (\ref{eq:vqdef}). There exist $R,C>0$ such that 
\begin{align*}
\left|\nabla_{\bar{y}}^{\tau}\left(v_{i}(y)-U(y)-\delta_{i}^{2}\gamma_{x_{i}}(y)\right)\right| & \le C\delta_{i}^{3}(1+|y|)^{5-\tau-n}\\
\left|y_{n}\frac{\partial}{\partial_{n}}\left(v_{i}(y)-U(y)-\delta_{i}^{2}\gamma_{x_{i}}(y)\right)\right| & \le C\delta_{i}^{3}(1+|y|)^{5-n}
\end{align*}
for $|y|\le\frac{R}{2\delta_{i}}$. Here $\tau=0,1,2$ and $\nabla_{\bar{y}}^{\tau}$
is the differential operator of order $\tau$ with respect the first
$n-1$ variables. 
\end{prop}

\begin{proof}
In analogy with Lemma \ref{lem:coreLemma}, we set 
\[
w_{i}(y):=v_{i}(y)-U(y)-\delta_{i}^{2}\gamma_{x_{i}}(y)\text{ for }|y|\le R/\delta_{i},
\]
and we have that $w_{i}$ satisfies (\ref{eq:wi}) where $b_{i}$
is defined as before, 
\begin{align*}
Q_{i}= & -\left\{ \left(L_{\hat{g}_{i}}-\Delta\right)(U+\delta_{i}^{2}\gamma_{x_{i}})+\delta_{i}^{2}\Delta\gamma_{x_{i}}\right\} ,\\
A_{i}= & Q_{i}+\varepsilon_{1,i}\delta_{i}^{2}\alpha_{i}(\delta_{i}y)v_{i}(y),\\
\bar{Q}_{i}= & -\left\{ (n-2)(U+\delta_{i}^{2}\gamma_{x_{i}})^{\frac{n}{n-2}}-(n-2)U^{\frac{n}{n-2}}-n\delta_{i}^{2}U^{\frac{2}{n-2}}\gamma_{x_{i}}-\frac{n-2}{2}h_{\hat{g}_{i}}(U+\delta_{i}^{2}\gamma_{x_{i}})\right\} ,\\
F_{i}= & \bar{Q}_{i}+\varepsilon_{2,i}\delta_{i}\beta_{i}(\delta_{i}y)v_{i}(y).
\end{align*}
As before, $b_{i}$ satisfies inequality (\ref{eq:b2}) and 
\begin{align}
A_{i} & =O(\delta_{i}^{3}\left(1+|y|\right)^{3-n})+O(\varepsilon_{1,i}\delta_{i}^{2}\left(1+|y|\right)^{2-n}),\label{eq:Q-1}\\
F_{i} & =O(\delta_{i}^{4}\left(1+|y|\right)^{5-n})+O(\varepsilon_{2,i}\delta_{i}\left(1+|y|\right)^{2-n}).\label{eq:Q-2}
\end{align}
We define the Green function $G_{i}$ as in the previous lemma and
again, by Green's formula, by (\ref{eq:Q-1}), (\ref{eq:Q-2}), and
Lemmas \ref{lem:coreLemma} and \ref{lem:taui}, we have
\begin{eqnarray}
|w_{i}(y)|\le C\delta_{i}^{3}\text{ on }B_{R/\delta_{i}}^{+} & \text{ and } & |w_{i}(\xi)|\le C\delta_{i}^{n-2}\text{ on }\partial^{+}B_{R/\delta_{i}}^{+}.\label{eq:w_i-improved}
\end{eqnarray}
By this we show that $\int_{\partial'B_{\frac{R}{\delta_{i}}}^{+}}|\bar{\xi}-y|^{2-n}b_{i}(\xi)w_{i}(\xi)d\bar{\xi})\le\delta_{i}^{3}(1+|y|)^{-1}$,
while one can manage the other terms in Green's formula as in the
previous lemmas. So we obtain
\begin{equation}
|w_{i}(y)|\le C\delta_{i}^{3}(1+|y|)^{-1}\text{ for }|y|\le\frac{R}{2\delta_{i}}.\label{eq:wi-improv2}
\end{equation}
Now we can iterate the procedure until we reach 
\begin{equation}
|w_{i}(y)|\le C\delta_{i}^{3}(1+|y|)^{5-n}\text{ for }|y|\le\frac{R}{2\delta_{i}},\label{eq:wi-improv2-1-1}
\end{equation}
which proves the first claim for $\tau=0$. The other claims follow
similarly.
\end{proof}

\subsection{Sign estimates of Pohozaev identity terms\label{sec:Sign-estimates}}

In this section, we want to estimate $P(u_{i},r)$, where $\left\{ u_{i}\right\} _{i}$
is a family of solutions of (\ref{eq:Prob-i}) which has an isolated
simple blow up point $x_{i}\rightarrow x_{0}$. This estimate, given
in the following Proposition \ref{prop:segno}, is a crucial point
for the proof of the vanishing of the Weyl tensor at an isolated simple
blow up point.

Since the leading term of $P(u_{i},r)$ will be $-\int_{B_{r/\delta_{i}}^{+}}\left(y^{b}\partial_{b}u+\frac{n-2}{2}u\right)\left[(L_{\hat{g}_{i}}-\Delta)v\right]dy$,
we set
\begin{equation}
R(u,v)=-\int_{B_{r/\delta_{i}}^{+}}\left(y^{b}\partial_{b}u+\frac{n-2}{2}u\right)\left[(L_{\hat{g}_{i}}-\Delta)v\right]dy,\label{eq:Ruv}
\end{equation}
and we recall the following result
\begin{lem}
\label{lem:R(U,U)}For $n\ge8$ we have 
\[
R(U+\delta^{2}\gamma_{q},U+\delta^{2}\gamma_{q})=\left\{ \begin{array}{c}
\delta^{4}\frac{(n-2)\omega_{n-2}I_{n}^{n}}{(n-1)(n-3)(n-5)(n-6)}\left[\frac{\left(n-2\right)}{6}|\bar{W}(q)|^{2}+\frac{4(n-8)}{(n-4)}R_{ninj}^{2}(q)\right]\\
-2\delta^{4}\int_{\mathbb{R}_{+}^{n}}\gamma_{q}\Delta\gamma_{q}dy+o(\delta^{4})\text{ for }n>8\\
\\
\delta_{i}^{4}\omega_{6}I_{8}^{8}\left[\frac{1}{35}|\bar{W}(q)|^{2}+\frac{1089}{34020}R_{8i8j}^{2}(q)\right]+o(\delta^{4})\text{ for }n=8
\end{array}\right..
\]
\end{lem}

\begin{proof}
For the proof we refer to \cite{GM20} for the case $n>8$ and to
\cite{GMsub} for the case $n=8$.
\end{proof}
\begin{prop}
\label{prop:segno}Let $x_{i}\rightarrow x_{0}$ be an isolated simple
blow-up point for $u_{i}$ solutions of (\ref{eq:Prob-i}). Let $\alpha,\beta<0$
and $n\ge8$. Then, fixed $r$, we have, for $i$ large 
\begin{align*}
\hat{P}(u_{i},r)\ge & \delta_{i}^{4}\left[C_{1}|\bar{W}(x_{i})|^{2}+C_{2}R_{nlnj}^{2}(x_{i})\right]+o(\delta_{i}^{4})
\end{align*}
where $C_{1},C_{2}>0$.
\end{prop}

\begin{proof}
We recall that the definition of $\hat{P}$ is given in Theorem \ref{thm:poho}
and we take $v_{i}(y)$ as in (\ref{eq:deltai}). By Proposition \ref{prop:stimawi}
and by (\ref{eq:gradvq}) of Lemma \ref{lem:vq}, for $|y|<R/\delta_{i}$
we have
\[
\left|v_{i}(y)-U(y)\right|=O(\delta_{i}^{3}(1+|y|^{5-n})+O(\delta_{i}^{2}(1+|y|^{4-n})=O(\delta_{i}^{2}(1+|y|^{4-n}),
\]
\[
\left|y_{a}\partial_{a}v_{i}(y)-y_{a}\partial_{a}U(y)\right|=O(\delta_{i}^{3}(1+|y|^{5-n})+O(\delta_{i}^{2}(1+|y|^{4-n})=O(\delta_{i}^{2}(1+|y|^{4-n}),
\]
so
\[
\int\limits _{B_{r}^{+}}\left(y^{a}\partial_{a}u_{i}+\frac{n-2}{2}u_{i}\right)\varepsilon_{1,i}\alpha_{i}u_{i}dy=\varepsilon_{1,i}\delta_{i}^{2}\int\limits _{B_{r/\delta_{i}}^{+}}\left(y^{a}\partial_{a}v_{i}+\frac{n-2}{2}v_{i}\right)\alpha_{i}(\delta_{i}y)v_{i}dy+\varepsilon_{1,i}\delta_{i}^{2}o(\delta_{i}^{2})
\]
and, recalling that $\alpha_{i}(\delta_{i}y)\rightarrow\alpha(x_{0})<0$
and proceeding as in Proposition \ref{prop:eps-i} we get
\[
\lim_{i\rightarrow\infty}\int\limits _{B_{r}^{+}}\left(y^{a}\partial_{a}v_{i}+\frac{n-2}{2}v_{i}\right)\alpha_{i}v_{i}dy=\frac{n-2}{2}\alpha(x_{0})\int\limits _{\mathbb{R}_{+}^{n}}\frac{1-|y|^{2}}{\left[(1+y_{n})^{2}+|\bar{y}|^{2}\right]^{n-1}}dy>0.
\]
Analogously 
\[
\int\limits _{\partial'B_{r}^{+}}\left(\bar{y}^{k}\partial_{k}u_{i}+\frac{n-2}{2}u_{i}\right)\varepsilon_{2,i}\beta_{i}u_{i}d\bar{y}=\varepsilon_{2,i}\delta_{i}\int\limits _{\partial'B_{r/\delta_{i}}^{+}}\left(\bar{y}^{k}\partial_{k}v_{i}+\frac{n-2}{2}v_{i}\right)\beta_{i}(\delta_{i}y)Ud\bar{y}+\varepsilon_{2,i}\delta_{i}O(\delta_{i}^{2})
\]
and again we get
\[
\lim_{i\rightarrow\infty}\int\limits _{\partial'B_{r/\delta_{i}}^{+}}\left(\bar{y}^{k}\partial_{k}v_{i}+\frac{n-2}{2}v_{i}\right)\beta_{i}(\delta_{i}y)v_{i}d\bar{y}=\frac{n-2}{2}\beta(x_{0})\int\limits _{\mathbb{R}^{n-1}}\frac{1-|\bar{y}|^{2}}{\left[1+|\bar{y}|^{2}\right]^{n-1}}d\bar{y}>0.
\]
Finally, since $h_{g_{i}}(\delta_{i}y)=O(\delta_{i}^{4}|y|^{4})$
we have$\int_{\partial'B_{r/\delta_{i}}^{+}}\left(y^{b}\partial_{b}v_{i}+\frac{n-2}{2}v_{i}\right)h_{g_{i}}(\delta_{i}y)v_{i}d\bar{y}=O(\delta_{i}^{5})$.
So, for $i$ sufficiently large we obtain 
\[
\hat{P}(u_{i},r)\ge-\int_{B_{r/\delta_{i}}^{+}}\left(y^{b}\partial_{b}v_{i}+\frac{n-2}{2}v_{i}\right)\left[(L_{\hat{g}_{i}}-\Delta)v_{i}\right]dy+O(\delta_{i}^{5}).
\]
Now define, in analogy with Proposition \ref{prop:stimawi}, 
\[
w_{i}(y):=v_{i}(y)-U(y)-\delta_{i}^{2}\gamma_{x_{i}}(y).
\]
Recalling (\ref{eq:Ruv}), we have 
\begin{align*}
\hat{P}(u_{i},r) & \ge R(U,U)+R(U,\delta_{i}^{2}\gamma_{x_{i}})+R(\delta_{i}^{2}\gamma_{x_{i}},U)+R(w_{i},U)+R(U,w_{i})\\
 & +R(w_{i,}w_{i})+R(\delta_{i}^{2}\gamma_{x_{i}},\delta_{i}^{2}\gamma_{x_{i}})+R(w_{i},\delta_{i}^{2}\gamma_{x_{i}})+R(\delta_{i}^{2}\gamma_{x_{i}},w_{i})+O(\delta_{i}^{5}).
\end{align*}
By \cite{GM20} we have that $\hat{P}(u_{i},r)\ge R(U,U)+R(U,\delta_{i}^{2}\gamma_{x_{i}})+R(\delta_{i}^{2}\gamma_{x_{i}},U)+R(\delta_{i}^{2}\gamma_{x_{i}},\delta_{i}^{2}\gamma_{x_{i}})+O(\delta_{i}^{5})$
and by Lemma \ref{lem:R(U,U)} we conclude the proof.
\end{proof}
\begin{prop}
\label{prop:7.1}Assume $n\ge8$ and $\alpha,\beta<0$. Let $x_{i}\rightarrow x_{0}$
be an isolated simple blow-up point for $u_{i}$ solutions of (\ref{eq:Prob-i}).
Then $|W(x_{0})|=0.$
\end{prop}

\begin{proof}
By Proposition \ref{prop:4.3} and Proposition \ref{prop:Lemma 4.4},
and since $M_{i}=\delta_{i}^{\frac{2-n}{2}}$ we have,
\begin{align*}
P(u_{i},r):= & \frac{1}{M_{i}^{2\lambda_{i}}}\int\limits _{\partial^{+}B_{r}^{+}}\left(\frac{n-2}{2}M_{i}^{\lambda_{i}}u_{i}\frac{\partial M_{i}^{\lambda_{i}}u_{i}}{\partial r}-\frac{r}{2}|\nabla M_{i}^{\lambda_{i}}u_{i}|^{2}+r\left|\frac{\partial M_{i}^{\lambda_{i}}u_{i}}{\partial r}\right|^{2}\right)d\sigma_{r}\\
 & +\frac{r(n-2)^{2}}{\left(n-1\right)M_{i}^{\lambda_{i}\frac{2(n-1)}{n-2}}}\int\limits _{\partial(\partial'B_{r}^{+})}\left(M_{i}^{\lambda_{i}}u_{i}\right)^{\frac{2(n-1)}{n-2}}d\bar{\sigma}_{g}.\\
\le & \frac{C}{M_{i}^{\lambda_{i}\frac{2(n-1)}{n-2}}}\le C\delta_{i}^{(n-1)\lambda_{i}}\le C\delta_{i}^{n-2}.
\end{align*}
On the other hand recalling Proposition \ref{prop:segno} and Theorem
\ref{thm:poho} we have 

\[
P(u_{i},r)=\hat{P}(u_{i},r)\ge\delta_{i}^{4}\left[C_{1}|\bar{W}(x_{i})|^{2}+C_{2}R_{nlnj}^{2}(x_{i})\right]+o(\delta_{i}^{4}),
\]
so we get $\left[C_{1}|\bar{W}(x_{i})|^{2}+C_{2}R_{nlnj}^{2}(x_{i})\right]\le\delta_{i}^{2}$.
Recalling that when the boundary is umbilic $W(q)=0$ if and only
if $\bar{W}(q)=0$ and $R_{nlnj}(q)=0$ (see \cite[page 1618]{M1})
we conclude the proof. 
\end{proof}
\begin{rem}
\label{rem:P'}Let $x_{i}\rightarrow x_{0}$ be an isolated blow up
point for $u_{i}$ solutions of (\ref{eq:Prob-i}). We set 
\begin{equation}
P'\left(u,r\right):=\int\limits _{\partial^{+}B_{r}^{+}}\left(\frac{n-2}{2}u\frac{\partial u}{\partial r}-\frac{r}{2}|\nabla u|^{2}+r\left|\frac{\partial u}{\partial r}\right|^{2}\right)d\sigma_{r},\label{eq:P'def}
\end{equation}
 so
\[
P(u_{i},r)=P'(u_{i},r)+\frac{r(n-2)^{2}}{(n-1)}\int\limits _{\partial(\partial'B_{r}^{+})}u_{i}^{\frac{2(n-1)}{n-2}}d\bar{\sigma}_{g}
\]
and, keeping in mind that for $i$ large $M_{i}u_{i}\le C|y|^{2-n}$
by Proposition \ref{prop:4.3}, we have 
\begin{equation}
\left|r\int\limits _{\partial(\partial'B_{r}^{+})}u_{i}^{\frac{2(n-1)}{n-2}}d\bar{\sigma}_{g}\right|\le\frac{Cr}{M_{i}^{\frac{2(n-1)}{n-2}}}\int_{\begin{array}{c}
y_{n}=0\\
|\bar{y}|=r
\end{array}}\frac{1}{|y|^{2(n-1)}}d\bar{\sigma}_{g}\le\frac{C(r)}{M_{i}^{\frac{2(n-1)}{n-2}}}=C(r)\delta_{i}^{n-2}\label{eq:5.15-1}
\end{equation}
for $i$ sufficiently large.

Using Proposition \ref{prop:segno}, (\ref{eq:5.15-1}), and since
$n\ge8$ we get 
\begin{equation}
P'(u_{i},r)=P(u_{i},r)-\frac{r(n-2)^{2}}{(n-1)}\int\limits _{\partial(\partial'B_{r}^{+})}u_{i}^{\frac{2(n-1)}{n-2}}d\bar{\sigma}_{g}\ge C\delta_{i}^{4}+o(\delta^{4})\label{eq:stimaP'}
\end{equation}
where $C>0$.
\end{rem}

\begin{prop}
\label{prop:isolato->semplice}Let $x_{i}\rightarrow x_{0}$ be an
isolated blow up point for $u_{i}$ solutions of (\ref{eq:Prob-i}).
Assume $n\ge8$ and $|W(x_{0})|\neq0$. Then $x_{0}$ is isolated
simple. 
\end{prop}

For the proof of this Lemma we refer to \cite{Al,GM20} 

\subsection{Proof of Theorem \ref{thm:main}\label{sec:Main-Proof}}

Before the proof of Theorem, we summarize a result which proves that
only isolated blow up points may occur to a blowing up sequence of
solution. For the proof of this result we refer to \cite[Proposition 5.1]{LZ},
\cite[Lemma 3.1]{SZ}, \cite[Proposition 1.1]{HL}, \cite[Proposition 4.2]{Al}
for the first claims, \cite{GM20} for the last claim when $n>8$
and to \cite{GMsub} in the case $n=8$.
\begin{prop}
\label{prop:4.2}Given $K>0$ and $R>0$ there exist two constants
$C_{0},C_{1}>0$ (depending on $K$, $R$ and $(M,g)$) such that
if $u$ is a solution of 
\begin{equation}
\left\{ \begin{array}{cc}
L_{g}u-\varepsilon_{1}\alpha=0 & \text{ in }M\\
B_{g}u-\varepsilon_{2}\beta u+(n-2)u^{\frac{n}{n-2}}=0 & \text{ on }\partial M
\end{array}\right.\label{eq:Prob-p}
\end{equation}
and $\max_{\partial M}u>C_{0}$, then there exist $q_{1},\dots,q_{N}\in\partial M$,
with $N=N(u)\ge1$ with the following properties: for $j=1,\dots,N$ 
\begin{enumerate}
\item set $r_{j}:=Ru(q_{j})^{1-p}$ then $\left\{ B_{r_{j}}\cap\partial M\right\} _{j}$
are a disjoint collection;
\item we have $\left|u(q_{j})^{-1}u(\psi_{j}(y))-U(u(q_{j})^{p-1}y)\right|_{C^{2}(B_{2r_{j}}^{+})}<K$
(here $\psi_{j}$ are the Fermi coordinates at point $q_{j}$;
\item we have
\begin{align*}
u(x)d_{\bar{g}}\left(x,\left\{ q_{1},\dots,q_{n}\right\} \right)^{\frac{1}{p-1}}\le C_{1} & \text{ for all }x\in\partial M\\
u(q_{j})d_{\bar{g}}\left(q_{j},q_{k}\right)^{\frac{1}{p-1}}\ge C_{0} & \text{ for any }j\neq k.
\end{align*}
\end{enumerate}
In addition, if $n\ge8$ and $W(x)\neq0$ for any $x\in\partial M$,
there exists $d=d(K,R)$ such that
\[
\min_{\begin{array}{c}
i\neq j\\
1\le i,j\le N(u)
\end{array}}d_{\bar{g}}(q_{i}(u),q_{j}(u))\ge d.
\]
Here $\bar{g}$ is the geodesic distance on $\partial M$.
\end{prop}

We prove now the main result
\begin{proof}[Proof of Theorem \ref{thm:main}]
. By contradiction, suppose that $x_{i}\rightarrow x_{0}$ is a blowup
point for $u_{i}$ solutions of (\ref{eq:Prob-2}). Let $q_{1}^{i},\dots q_{N(u_{i})}^{i}$
the sequence of points given by Proposition \ref{prop:4.2}. By Claim
3 of Proposition \ref{prop:4.2} there exists a sequence of indices
$k_{i}\in1,\dots N$ such that $d_{\bar{g}}\left(x_{i},q_{k_{i}}^{i}\right)\rightarrow0$.
Up to relabeling, we say $k_{i}=1$ for all $i$. Then also $q_{1}^{i}\rightarrow x_{0}$
is a blow up point for $u_{i}$. By Proposition \ref{prop:4.2} and
Proposition \ref{prop:isolato->semplice} we have that $q_{1}^{i}\rightarrow x_{0}$
is an isolated simple blow up point for $u_{i}$. Then by Proposition
\ref{prop:7.1} we deduce that $W(x_{0})=0$, contradicting the assumption
of the theorem. This concludes the proof.
\end{proof}

\section{The non compactness result\label{sec:The-non-compactness} }

In this section we perform the Ljapunov-Schmidt finite dimensional
reduction, which relies on three steps. First, we start finding a
solution of the infinite dimensional problem (\ref{eq:Pibot}) with
the ansatz $\Lambda_{q}u=\tilde{W}_{\delta,q}+\delta^{2}\tilde{V}_{\delta,q}+\tilde{\phi}$
where $\tilde{\phi}\in\tilde{K}_{\delta,q}^{\bot}$. This is done
in subsection \ref{subsec:The-finite-dimensional}. Then, we study
the finite dimensional reduced problem in subsection \ref{subsec:The-reduced-functional},
and in the last subsection we give the proof of Theorem \ref{thm:main2}.

\subsection{\label{subsec:The-finite-dimensional}The finite dimensional reduction}

Let us define the linear operator $L:\tilde{K}_{\delta,q}^{\bot}\rightarrow\tilde{K}_{\delta,q}^{\bot}$
as
\begin{equation}
L(\tilde{\phi})=\tilde{\Pi}^{\bot}\left\{ \tilde{\phi}-i_{\alpha}^{*}\left(f'(\tilde{W}_{\delta,q}+\delta^{2}\tilde{V}_{\delta,q})[\tilde{\phi}]\right)\right\} ,\label{eq:defL}
\end{equation}
and let us define a nonlinear term $N(\tilde{\phi})$ and a remainder
term R as 
\begin{align}
N(\tilde{\phi})= & \tilde{\Pi}^{\bot}\left\{ i_{\alpha}^{*}\left(f(\tilde{W}_{\delta,q}+\delta^{2}\tilde{V}_{\delta,q}+\tilde{\phi})-f(\tilde{W}_{\delta,q}+\delta^{2}\tilde{V}_{\delta,q})-f'(\tilde{W}_{\delta,q}+\delta^{2}\tilde{V}_{\delta,q})[\tilde{\phi}]\right)\right\} \label{eq:defN}\\
R= & \tilde{\Pi}^{\bot}\left\{ i_{\alpha}^{*}\left(f(\tilde{W}_{\delta,q}+\delta^{2}\tilde{V}_{\delta,q})\right)-\tilde{W}_{\delta,q}-\delta^{2}\tilde{V}_{\delta,q}\right\} ,\label{eq:defR}
\end{align}
With these operators the infinite dimensional equation (\ref{eq:Pibot})
becomes

\[
L(\tilde{\phi})=N(\tilde{\phi})+R-\tilde{\Pi}^{\bot}\left\{ i_{\alpha}^{*}\left(\varepsilon_{2}\beta(\tilde{W}_{\delta,q}+\delta^{2}\tilde{V}_{\delta,q}+\tilde{\phi})\right)\right\} .
\]
In this subsection we will find, for any $\delta,q$ given, a function
$\tilde{\phi}$ which solves equation (\ref{eq:Pibot}). 
\begin{lem}
\label{lem:R}It holds 
\[
\|R\|_{g}=\left\{ \begin{array}{cc}
O\left(\delta^{3}\log\delta\right)+O(\varepsilon_{1}\delta^{2}) & \text{ if }n=8\\
O\left(\delta^{3}\right)+O(\varepsilon_{1}\delta^{2}) & \text{\text{ if }}n>8
\end{array}\right..
\]
\end{lem}

\begin{proof}
Several estimates for this proof has been calculated in \cite{GMwip},
which we refer to. We report here only the main steps.

Take the unique $\Gamma$ such that 
\[
\Gamma=i_{\alpha}^{*}\left(f(\tilde{W}_{\delta,q}+\delta^{2}\tilde{V}_{\delta,q})\right),
\]
that is the unique $\Gamma$ which solves 
\[
\left\{ \begin{array}{ccc}
-\Delta_{g}\Gamma+\frac{n-2}{4(n-1)}R_{g}\Gamma+\varepsilon_{1}\alpha\Gamma=0 &  & \text{on }M;\\
\frac{\partial\Gamma}{\partial\nu}+\frac{n-2}{2}h_{g}\Gamma=(n-2)\left((\tilde{W}_{\delta,q}+\delta^{2}\tilde{V}_{\delta,q})^{+}\right)^{\frac{n}{n-2}} &  & \text{on \ensuremath{\partial}}M.
\end{array}\right.
\]
Let us call $a:=\frac{n-2}{4(n-1)}R_{g}$. We have, by (\ref{eq:prodscal})
that 
\begin{align*}
\|R\|_{g}^{2} & \le\left\Vert i_{g}^{*}\left(f(\tilde{W}_{\delta,q}+\delta^{2}\tilde{V}_{\delta,q}\right)-\tilde{W}_{\delta,q}-\delta^{2}\tilde{V}_{\delta,q}\right\Vert _{g}^{2}=\|\Gamma-\tilde{W}_{\delta,q}-\delta^{2}\tilde{V}_{\delta,q}\|_{g}^{2}\\
 & =\int_{M}\left[\Delta_{g}(\tilde{W}_{\delta,q}+\delta^{2}\tilde{V}_{\delta,q})-a(\tilde{W}_{\delta,q}+\delta^{2}\tilde{V}_{\delta,q})\right](\Gamma-\tilde{W}_{\delta,q}-\delta^{2}\tilde{V}_{\delta,q})d\mu_{g}\\
 & -\int_{M}\varepsilon_{1}\alpha(\tilde{W}_{\delta,q}+\delta^{2}\tilde{V}_{\delta,q})(\Gamma-\tilde{W}_{\delta,q}-\delta^{2}\tilde{V}_{\delta,q})d\mu_{g}\\
 & -\int_{\partial M}h_{g}(\tilde{W}_{\delta,q}+\delta^{2}\tilde{V}_{\delta,q})(\Gamma-\tilde{W}_{\delta,q}-\delta^{2}\tilde{V}_{\delta,q})d\sigma_{g}\\
 & +\int_{\partial M}\left[f(\tilde{W}_{\delta,q}+\delta^{2}\tilde{V}_{\delta,q})-\frac{\partial}{\partial\nu}(\tilde{W}_{\delta,q}+\delta^{2}\tilde{V}_{\delta,q})\right](\Gamma-\tilde{W}_{\delta,q}-\delta^{2}\tilde{V}_{\delta,q})d\sigma_{g}\\
 & =:I_{1}+I_{2}+I_{3}+I_{4}.
\end{align*}
We have that 
\begin{align*}
I_{3} & =\int_{\partial M}h_{\tilde{g}}(W_{\delta,q}+\delta^{2}V_{\delta,q})(\Lambda_{q}^{-1}R)d\sigma_{\tilde{g}}\\
 & \le C|h_{\tilde{g}}(W_{\delta,q}+\delta^{2}V_{\delta,q})|_{L_{\tilde{g}}^{\frac{2(n-1)}{n}}(\partial M)}\|\Lambda_{q}^{-1}R\|_{\tilde{g}},
\end{align*}
and, by change of variables and by (\ref{eq:hij}), we get
\[
|h_{\tilde{g}}(W_{\delta,q}+\delta^{2}V_{\delta,q})|_{L_{\tilde{g}}^{\frac{2(n-1)}{n}}(\partial M)}=\left\{ \begin{array}{cc}
O(\delta^{3}\log\delta) & \text{ if }n=8\\
O(\delta^{3}) & \text{ if }n>8
\end{array}\right..
\]
Similarly for $I_{1}$ we have 
\[
I_{1}\le\left|\Delta_{\tilde{g}}(W_{\delta,q}+\delta^{2}V_{\delta,q})-\tilde{a}(W_{\delta,q}+\delta^{2}V_{\delta,q})\right|_{L_{\tilde{g}}^{\frac{2n}{n+2}}(M)}\|\Lambda_{q}^{-1}R\|_{\tilde{g}}
\]
 Since $R_{\tilde{g}}(0)=0$ (see \cite[page 1609]{M1}), we get 
\[
\left|\tilde{a}(W_{\delta,q}+\delta^{2}V_{\delta,q})\right|_{L_{\tilde{g}}^{\frac{2n}{n+2}}(M)}=\left\{ \begin{array}{cc}
O(\delta^{3}\log\delta) & \text{ if }n=8\\
O(\delta^{3}) & \text{ if }n>8
\end{array}\right.
\]
and, using the expansion of the metric $\tilde{g}$ and (\ref{eq:vqdef}),
one can show that 
\[
\left|\Delta_{\tilde{g}}(W_{\delta,q}+\delta^{2}V_{\delta,q})\right|_{L_{\tilde{g}}^{\frac{2n}{n+2}}(M)}=\left\{ \begin{array}{cc}
O(\delta^{3}\log\delta) & \text{ if }n=8\\
O(\delta^{3}) & \text{ if }n>8
\end{array}\right.;
\]
thus we get 
\[
I_{1}+I_{3}=\left\{ \begin{array}{cc}
O(\delta^{3}\log\delta)\|R\|_{g} & \text{ if }n=8\\
O(\delta^{3})\|R\|_{g} & \text{ if }n>8
\end{array}\right..
\]
For the integral $I_{4}$ we have
\begin{align*}
I_{4}\le & C(n-2)\left|\left((W_{\delta,q}+\delta^{2}V_{\delta,q})^{+}\right)^{\frac{n}{n-2}}-\left(W_{\delta,q}\right)^{\frac{n}{n-2}}-\delta^{2}\frac{\partial}{\partial\nu}V_{\delta,q}\right|_{L_{\tilde{g}}^{\frac{2(n-1)}{n}}(\partial M)}\left\Vert R\right\Vert _{g}\\
 & +C\left|(n-2)\left(W_{\delta,q}\right)^{\frac{n}{n-2}}-\frac{\partial}{\partial\nu}W_{\delta,q}\right|_{L_{\tilde{g}}^{\frac{2(n-1)}{n}}(\partial M)}\left\Vert R\right\Vert _{g}.
\end{align*}
and, since $U$ solves (\ref{eq:Udelta}), we get immediately
\begin{equation}
\left|(n-2)\left(W_{\delta,q}\right)^{\frac{n}{n-2}}-\frac{\partial}{\partial\nu}W_{\delta,q}\right|_{L_{\tilde{g}}^{\frac{2(n-1)}{n}}(\partial M)}=O(\delta^{3}).\label{eq:bordo2}
\end{equation}
Estimating the other terms requires more care, but, expanding $\left((U+\delta^{2}\gamma_{q})^{+}\right)^{\frac{n}{n-2}}$
near $U$, using  (\ref{eq:vqdef}) and the decay estimates (\ref{eq:gradvq}),
one can show that (see \cite{GMwip} for all the details )
\begin{multline}
\left|\left((W_{\delta,q}+\delta^{2}V_{\delta,q})^{+}\right)^{\frac{n}{n-2}}-\left(W_{\delta,q}\right)^{\frac{n}{n-2}}-\delta^{2}\frac{\partial}{\partial\nu}V_{\delta,q}\right|_{L_{\tilde{g}}^{\frac{2(n-1)}{n}}(\partial M)}=O(\delta^{3}).\label{eq:bordo3}
\end{multline}
Thus by (\ref{eq:bordo2}) and (\ref{eq:bordo3}) we have 
\[
I_{4}=O(\delta^{3}).
\]
For $I_{2}$ we have 
\[
I_{2}\le\varepsilon_{1}\left|\tilde{a}(W_{\delta,q}+\delta^{2}V_{\delta,q})\right|_{L_{\tilde{g}}^{\frac{2n}{n+2}}(M)}\|\Lambda_{q}^{-1}R\|_{\tilde{g}}
\]
Now by change of variables we have 
\[
\left|W_{\delta,q}+\delta^{2}V_{\delta,q}\right|_{L_{\tilde{g}}^{\frac{2n}{n+2}}(M)}=O(\delta^{2})\left|(1+|x|)^{2-n}\right|_{L^{\frac{2n}{n+2}}(B_{1/\delta})}=O(\delta^{2})
\]
 so 
\[
I_{2}=O(\varepsilon_{1}\delta^{2})\left\Vert R\right\Vert _{g}
\]
 which completes the proof.
\end{proof}
The following lemma is a standard tool in finite dimensional reduction,
so we refer to \cite{EPV14,MP09} for its proof.
\begin{lem}
\label{lem:L}Given $(\varepsilon_{1},\varepsilon_{2})$, for any
pair $(\delta,q)$ there exists a positive constant $C=C(\delta,q)$
such that for any $\varphi\in\tilde{K}_{\delta,q}^{\bot}$ it holds
\[
\|L(\varphi)\|_{g}\ge C\|\varphi\|_{g}.
\]
\end{lem}

It is also standard to prove that $N$ is a contraction, that is there
exists $\eta<1$ such that, for any $\varphi_{1},\varphi_{2}\in\tilde{K}_{\delta,q}^{\bot}$
it holds 
\begin{equation}
\|N(\varphi)\|_{g}\le\eta\|\varphi\|_{g}\text{ and }\|N(\varphi_{1})-N(\varphi_{2})\|_{g}\le\eta\|\varphi_{1}-\varphi_{2}\|_{g}\label{eq:Ncontr}
\end{equation}
By Lemma \ref{lem:R}, Lemma \ref{lem:L}, and (\ref{eq:Ncontr})
we prove the last result of this subsection.
\begin{prop}
\label{prop:EsistenzaPhi}Given $(\varepsilon_{1},\varepsilon_{2})$,
for any pair $(\delta,q)$ there exists a unique $\tilde{\phi}=\tilde{\phi}_{\delta,q}\in\tilde{K}_{\delta,q}^{\bot}$
which solves (\ref{eq:Pibot}) such that 
\[
\|\tilde{\phi}\|_{g}=\left\{ \begin{array}{cc}
O\left(\delta^{3}\log\delta+\varepsilon_{1}\delta^{2}+\varepsilon_{2}\delta\right) & \text{ if }n=8\\
O\left(\delta^{3}+\varepsilon_{1}\delta^{2}+\varepsilon_{2}\delta\right) & \text{\text{ if }}n>8
\end{array}\right..
\]
The map $q\mapsto\phi$ is $C^{1}$.
\end{prop}

\begin{proof}
By Lemma \ref{lem:L}, by (\ref{eq:Ncontr}) and by the properties
of $i_{\alpha}$, there exists $C>0$ such that 
\begin{multline*}
\left\Vert L^{-1}\left(N(\tilde{\phi})+R-\Pi^{\bot}\left\{ i_{\alpha}^{*}\left(\varepsilon_{2}\beta(\tilde{W}_{\delta,q}+\delta^{2}\tilde{V}_{\delta,q}+\tilde{\phi})\right)\right\} \right)\right\Vert _{g}\\
\le C\left((\eta\|\tilde{\phi}\|_{g}+\|R\|_{g}+\left\Vert i_{\alpha}^{*}\left(\varepsilon_{2}\beta(\tilde{W}_{\delta,q}+\delta^{2}\tilde{V}_{\delta,q}+\tilde{\phi})\right)\right\Vert _{g}\right).
\end{multline*}
Now it is easy to estimate that
\begin{align}
\left\Vert i_{\alpha}^{*}\left(\varepsilon_{2}\beta(\tilde{W}_{\delta,q}+\delta^{2}\tilde{V}_{\delta,q}+\tilde{\phi})\right)\right\Vert _{g} & \le\varepsilon_{2}\left(\left\Vert \tilde{W}_{\delta,q}+\delta^{2}\tilde{V}_{\delta,q}\right\Vert _{L_{g}^{\frac{2(n-1)}{n}}(\partial M)}+\left\Vert \tilde{\phi}\right\Vert _{g}\right)\nonumber \\
 & \le C\left(\varepsilon_{2}\delta+\varepsilon_{2}\left\Vert \tilde{\phi}\right\Vert _{g}\right).\label{eq:e2beta}
\end{align}
If $n>8$, by Lemma \ref{lem:R} and by the previous estimates, for
the map 
\[
T(\tilde{\phi}):=L^{-1}\left(N(\tilde{\phi})+R-\Pi^{\bot}\left\{ i_{\alpha}^{*}\left(\varepsilon_{2}\beta(\tilde{W}_{\delta,q}+\delta^{2}\tilde{V}_{\delta,q}+\tilde{\phi})\right)\right\} \right)
\]
it holds
\[
\|T(\tilde{\phi})\|_{g}\le C\left((\eta+\varepsilon_{2})\|\tilde{\phi}\|_{g}+\varepsilon_{2}\delta+\varepsilon_{1}\delta^{2}+\delta^{3}\right).
\]
It is possible to choose $\rho>0$ such that $T$ is a contraction
from the ball $\|\tilde{\phi}\|_{g}\le\rho(\varepsilon_{2}\delta+\varepsilon_{1}\delta^{2}+\delta^{3})$
in itself, so, by the fixed point Theorem, there exists a unique $\tilde{\phi}$
with $\|\tilde{\phi}\|_{g}=O(\varepsilon_{2}\delta+\varepsilon_{1}\delta^{2}+\delta^{3})$
which solves (\ref{eq:Pibot}). In addition by the implicit function
Theorem it is possible to prove the regularity of the map $q\mapsto\tilde{\phi}$.
The case $n=8$ follows verbatim.
\end{proof}

\subsection{\label{subsec:The-reduced-functional}The reduced functional}

Once we solved (\ref{eq:Pibot}), we show that we can find a critical
point of $J_{g}\left(\tilde{W}_{\delta,q}+\delta^{2}\tilde{V}_{\delta,q}+\tilde{\phi}\right)$
by solving a finite dimensional problem depending only on $(\delta,q)$.
\begin{lem}
\label{lem:JWpiuPhi}Assume $n\ge8$. It holds 
\begin{multline*}
\left|J_{g}\left(\tilde{W}_{\delta,q}+\delta^{2}\tilde{V}_{\delta,q}+\tilde{\phi}\right)-J_{g}\left(\tilde{W}_{\delta,q}+\delta^{2}\tilde{V}_{\delta,q}\right)\right|\\
=O\left(\left\Vert \tilde{\phi}\right\Vert _{g}^{2}+\delta^{2}\left\Vert \tilde{\phi}\right\Vert _{g}+\varepsilon_{2}\delta\left\Vert \tilde{\phi}\right\Vert _{g}\right)
\end{multline*}
$C^{0}$-uniformly for $q\in\partial M$. 
\end{lem}

\begin{proof}
We have, for some $\theta\in(0,1)$ 
\begin{multline*}
\tilde{J}_{\tilde{g}}(W_{\delta,q}+\delta^{2}V_{\delta,q}+\phi)-\tilde{J}_{\tilde{g}}(W_{\delta,q}+\delta^{2}V_{\delta,q})=\tilde{J}_{\tilde{g}}'(W_{\delta,q}+\delta^{2}V_{\delta,q})[\phi]\\
+\frac{1}{2}\tilde{J}_{\tilde{g}}''(W_{\delta,q}+\delta^{2}V_{\delta,q}+\theta\phi)[\phi,\phi]\\
=\int_{M}\left(\nabla_{\tilde{g}}W_{\delta,q}+\delta^{2}\nabla_{\tilde{g}}V_{\delta,q}\right)\nabla_{\tilde{g}}\phi+\left(\frac{n-2}{4(n-1)}R_{\tilde{g}}+\varepsilon_{1}\tilde{\alpha}\right)\left(W_{\delta,q}+\delta^{2}V_{\delta,q}\right)\phi d\mu_{\tilde{g}}\\
-(n-2)\int_{\partial M}\left(\left(W_{\delta,q}+\delta^{2}V_{\delta,q}\right)^{+}\right)^{\frac{n}{n-2}}\phi d\sigma_{\tilde{g}}+\frac{n-2}{2}\int_{\partial M}h_{\tilde{g}_{q}}\left(W_{\delta,q}+\delta^{2}V_{\delta,q}\right)\phi d\sigma_{\tilde{g}}\\
+\int_{\partial M}\varepsilon_{2}\tilde{\beta}\left(W_{\delta,q}+\delta^{2}V_{\delta,q}\right)\phi d\sigma_{\tilde{g}}+\frac{1}{2}\|\phi\|_{\tilde{g}}^{2}\\
-\frac{n}{2}\int_{\partial M}\left(\left(W_{\delta,q}+\delta^{2}V_{\delta,q}+\theta\phi_{\delta,q}\right)^{+}\right)^{\frac{2}{n-2}}\phi_{\delta,q}^{2}d\sigma_{\tilde{g}_{q}}+\frac{1}{2}\int_{\partial M}\varepsilon_{2}\tilde{\beta}|\phi|^{2}d\sigma_{\tilde{g}}.
\end{multline*}
Now, by Holder inequality we have 
\[
\left|\int_{M}W_{\delta,q}\phi d\mu_{\tilde{g}}\right|\le C|W_{\delta,q}|_{L_{\tilde{g}}^{\frac{2n}{n+2}}}|\phi|_{L_{\tilde{g}}^{\frac{2n}{n-2}}}\le C\delta^{2}\|\phi\|_{\tilde{g}}
\]
and 
\[
\delta^{2}\left|\int_{M}V_{\delta,q}\phi d\mu_{\tilde{g}}\right|\le C\delta^{2}|V_{\delta,q}|_{L_{\tilde{g}}^{2}}|\phi|_{L_{\tilde{g}}^{2}}\le C\delta^{2}\|\phi\|_{\tilde{g}}.
\]
Immediately we have $\left|\int_{\partial M}\varepsilon_{2}\tilde{\beta}|\phi|^{2}d\sigma_{\tilde{g}}\right|\le C\varepsilon_{2}\|\phi\|_{\tilde{g}}^{2}$,
and following the proof of \cite[Lemma 8]{GMwip} we obtain that 
\begin{multline*}
\left|\int_{M}\left(\nabla_{\tilde{g}}W_{\delta,q}+\delta^{2}\nabla_{\tilde{g}}V_{\delta,q}\right)\nabla_{\tilde{g}}\phi-(n-2)\int_{\partial M}\left(\left(W_{\delta,q}+\delta^{2}V_{\delta,q}\right)^{+}\right)^{\frac{n}{n-2}}\phi d\sigma_{\tilde{g}}\right|\\
\le C\delta^{2}\|\phi\|_{\tilde{g}}
\end{multline*}
and 
\[
\left|\int_{\partial M}\left(\left(W_{\delta,q}+\delta^{2}V_{\delta,q}+\theta\phi\right)^{+}\right)^{\frac{2}{n-2}}\phi_{\delta,q}^{2}d\sigma\right|\le C\|\phi\|_{\tilde{g}}^{2}.
\]
Finally by (\ref{eq:hij}) we have 
\[
\left|\int_{\partial M}h_{\tilde{g}_{q}}\left(W_{\delta,q}+\delta^{2}V_{\delta,q}\right)\phi d\sigma_{\tilde{g}}\right|\le C\delta^{3}\|\phi\|_{\tilde{g}}
\]
 and, proceeding similarly to (\ref{eq:e2beta}), 
\[
\left|\int_{\partial M}\varepsilon_{2}\tilde{\beta}\left(W_{\delta,q}+\delta^{2}V_{\delta,q}\right)\phi d\sigma_{\tilde{g}}\right|\le C\varepsilon_{2}\delta\|\phi\|_{\tilde{g}}.
\]
This proves the first estimate. Using the result of Proposition \ref{prop:EsistenzaPhi}
we complete the proof.
\end{proof}
\begin{lem}
\label{lem:espansione}Let $n\ge8$. It holds
\[
J_{g}(\tilde{W}_{\delta,q}+\delta^{2}\tilde{V}_{\delta,q})=A+\varepsilon_{1}\delta^{2}\alpha(q)B+\varepsilon_{2}\delta\beta(q)C+\delta^{4}\varphi(q)+O(\varepsilon_{1}\delta^{4})+O(\varepsilon_{2}\delta^{3})+O(\delta^{5})
\]
where 
\begin{align*}
A= & \frac{1}{2}\int_{\mathbb{R}_{+}^{n}}|\nabla U(y)|^{2}dy-\frac{(n-2)^{2}}{2(n-1)}\int_{\mathbb{R}^{n-1}}U(\bar{y},0)^{\frac{2(n-1)}{(n-2)}}d\bar{y}\\
B= & \frac{1}{2}\alpha(q)\int_{\mathbb{R}^{n}}U(y)^{2}dy\\
C= & \frac{1}{2}\beta(q)\int_{\mathbb{R}^{n-1}}U(\bar{y},0)^{2}d\bar{y}\\
\varphi(q)= & \frac{1}{2}\int_{\mathbb{R}_{+}^{n}}v_{q}\Delta v_{q}dy-\frac{n-2}{96(n-1)}|\bar{W}(q)|^{2}\int_{\mathbb{R}_{+}^{n}}|\bar{y}|^{2}U^{2}(\bar{y},y_{n})dy\\
 & -\frac{(n-2)(n-8)}{2(n^{2}-1)}R_{ninj}^{2}(q)\int_{\mathbb{R}_{+}^{n}}\frac{y_{n}^{2}|\bar{y}|^{4}}{\left((1+y_{n})^{2}+|\bar{y}|^{2}\right)^{n}}dy.
\end{align*}
Here $\bar{W}(q)$ is the Weyl tensor restricted to boundary. 
\end{lem}

\begin{proof}
The main estimates of this proof are proved in \cite[Lemma 8]{GMP},
which we refer to for a detailed proof, here we limit ourselves to
estimate the perturbation terms. We have
\begin{align*}
J_{g}(\tilde{W}_{\delta,q}+\delta^{2}\tilde{V}_{\delta,q})= & \frac{1}{2}\int_{M}|\nabla_{\tilde{g}_{q}}(W_{\delta,q}+\delta^{2}V_{\delta,q})|^{2}d\mu_{\tilde{g}}+\frac{n-2}{8(n-1)}\int_{M}R_{\tilde{g}}(W_{\delta,q}+\delta^{2}V_{\delta,q})^{2}d\mu_{\tilde{g}}\\
 & +\frac{1}{2}\varepsilon_{1}\int_{M}\Lambda_{q}^{-\frac{4}{n-2}}\alpha(W_{\delta,q}+\delta^{2}V_{\delta,q})^{2}d\mu_{\tilde{g}}\\
 & +\frac{1}{2}\varepsilon_{2}\int_{\partial M}\Lambda_{q}^{-\frac{2}{n-2}}\beta(W_{\delta,q}+\delta^{2}V_{\delta,q})^{2}d\sigma_{\tilde{g}}\\
 & -\frac{(n-2)^{2}}{2(n-1)}\int_{\partial M}\left(W_{\delta,q}+\delta^{2}V_{\delta,q}\right)^{\frac{2(n-1)}{n-2}}d\sigma_{\tilde{g}_{q}}\\
 & +\frac{n-2}{4}\int_{\partial M}h_{\tilde{g}_{q}}(W_{\delta,q}+\delta^{2}V_{\delta,q})^{2}d\sigma_{\tilde{g}_{q}}
\end{align*}
We easily compute the terms involving $\varepsilon_{1}$ and $\varepsilon_{2}$,
taking in account that $\Lambda_{q}(q)=1$ and the expansion of the
volume form given by (\ref{eq:|g|}), getting
\[
\frac{1}{2}\varepsilon_{1}\int_{M}\Lambda_{q}^{-\frac{4}{n-2}}\alpha(W_{\delta,q}+\delta^{2}V_{\delta,q})^{2}d\mu_{\tilde{g}}=\frac{1}{2}\varepsilon_{1}\delta^{2}\alpha(q)\int_{\mathbb{R}^{n}}U(y)^{2}dy+O(\varepsilon_{1}\delta^{4}),
\]
and 
\[
\frac{1}{2}\varepsilon_{2}\int_{\partial M}\Lambda_{q}^{-\frac{2}{n-2}}\beta(W_{\delta,q}+\delta^{2}V_{\delta,q})^{2}d\mu_{\tilde{g}}=\frac{1}{2}\varepsilon_{2}\delta\beta(q)\int_{\mathbb{R}^{n-1}}U(\bar{y},0)^{2}d\bar{y}+O(\varepsilon_{2}\delta^{3}).
\]
The remaining terms are estimated in \cite[Lemma 8]{GMP}, and it
holds
\begin{multline*}
\frac{1}{2}\int_{M}|\nabla_{\tilde{g}_{q}}(W_{\delta,q}+\delta^{2}V_{\delta,q})|^{2}d\mu_{\tilde{g}}+\frac{n-2}{8(n-1)}\int_{M}R_{\tilde{g}}(W_{\delta,q}+\delta^{2}V_{\delta,q})^{2}d\mu_{\tilde{g}}\\
-\frac{(n-2)^{2}}{2(n-1)}\int_{\partial M}\left(W_{\delta,q}+\delta^{2}V_{\delta,q}\right)^{\frac{2(n-1)}{n-2}}d\sigma_{\tilde{g}_{q}}+\frac{n-2}{4}\int_{\partial M}h_{\tilde{g}_{q}}(W_{\delta,q}+\delta^{2}V_{\delta,q})^{2}d\sigma_{\tilde{g}_{q}}\\
=A+\delta^{4}\varphi(q)+O(\delta^{5})
\end{multline*}
which completes the proof. 
\end{proof}

\subsection{Proof of Theorem \ref{thm:main2}}

At first we provide a sign estimate for function $\varphi(q)$ defined
in the previous paragraph.
\begin{lem}
\label{lem:phi<0}Assume $n\ge8$ and that the Weyl tensor $W_{g}$
is not vanishing on $\partial M$. Then the function $\varphi(q)$
defined in Lemma \ref{lem:espansione} is strictly negative on $\partial M$. 
\end{lem}

\begin{proof}
We can write the function $\varphi(q)$ defined in Lemma \ref{lem:espansione}
as
\[
\varphi(q)=\frac{1}{2}\int_{\mathbb{R}_{+}^{n}}\gamma_{q}\Delta\gamma_{q}dtdz-C_{1}|\bar{W}(q)|^{2}-(n-8)C_{2}R_{ninj}^{2}(q),
\]
where $C_{1},C_{2}$ are positive constants. If $n>8$, since in umbilic
boundary manifolds $W(q)=0$ if and only if $\bar{W}(q)$ and $R_{ninj}(q)$
are both zero (see \cite[page 1618]{M1}), by our assumption at least
one among $\bar{|W}(q)|$ and $R_{nlnj}^{2}(q)$ is strictly positive.
Since by (\ref{new}) the term involving $\gamma_{q}$ is non positive,
the lemma is proved.

When $n=8$ the term involving $R_{8i8j}^{2}$ vanishes. However in
\cite{GMsub} a refined analysis of the term $\int_{\mathbb{R}_{+}^{n}}\gamma_{q}\Delta\gamma_{q}dtdz$
was performed, leading to the following improvement of estimate (\ref{new}):
\[
\int_{\mathbb{R}_{+}^{8}}\gamma_{q}\Delta\gamma_{q}dy\le-C_{3}R_{8i8j}^{2}(q),
\]
where $C_{3}>0$. This was possible by a more precise description
of function $\gamma_{q}$ as sum of an harmonic function with explicit
rational functions, proved in Lemma 19 of the cited paper.

Thus for $n=8$ we have 
\[
\varphi(q)\le-C_{1}|\bar{W}(q)|^{2}-C_{3}R_{8i8j}^{2}(q)<0,
\]
and the proof is complete.
\end{proof}
\begin{proof}[Proof of Theorem \ref{thm:main2}]
 We give a detailed proof in the case $\alpha>0$. The case $\beta>0$
is analogous and we will emphasize the difference at the end of the
proof. 

If $\alpha>0$ we choose 
\begin{align*}
\delta & =\sqrt{\lambda\varepsilon_{1}}\\
\varepsilon_{2} & =o(\varepsilon_{1}^{2})
\end{align*}
where $\lambda\in\mathbb{R}^{+}$. With this choice, by Lemma \ref{lem:JWpiuPhi},
we have that 
\[
\left|J_{g}\left(\tilde{W}_{\sqrt{\lambda\varepsilon_{1}},q}+\lambda\varepsilon_{1}\tilde{V}_{\sqrt{\lambda\varepsilon_{1}},q}+\tilde{\phi}\right)-J_{g}\left(\tilde{W}_{\sqrt{\lambda\varepsilon_{1}},q}+\lambda\varepsilon_{1}\tilde{V}_{\sqrt{\lambda\varepsilon_{1}},q}\right)\right|=o(\varepsilon_{1}^{2})
\]
and that, by Lemma \ref{lem:espansione}, 
\[
J_{g}\left(\tilde{W}_{\sqrt{\lambda\varepsilon_{1}},q}+\lambda\varepsilon_{1}\tilde{V}_{\sqrt{\lambda\varepsilon_{1}},q}\right)=A+\varepsilon_{1}^{2}\left(\lambda\alpha(q)B+\lambda^{2}\varphi(q)\right)+o(\varepsilon_{1}^{2}).
\]

We recall a result which is a key tool in Ljapunov-Schmidt procedure,
and which is proved, for instance, in \cite[Lemma 9]{GMP} and which
relies on the estimates of Lemma \ref{lem:JWpiuPhi}.
\begin{rem*}
Given $(\varepsilon_{1},\varepsilon_{2}),$ if $(\bar{\lambda},\bar{q})\in(0,+\infty)\times\partial M$
is a critical point for the reduced functional $I_{\varepsilon_{1},\varepsilon_{2}}(\lambda,q):=J_{g}\left(\tilde{W}_{\sqrt{\lambda\varepsilon_{1}},q}+\lambda\varepsilon_{1}\tilde{V}_{\sqrt{\lambda\varepsilon_{1}},q}+\tilde{\phi}\right)$,
then the function $\tilde{W}_{\sqrt{\bar{\lambda}\varepsilon_{1}},\bar{q}}+\bar{\lambda}\varepsilon_{1}\tilde{V}_{\sqrt{\lambda\varepsilon_{1}},\bar{q}}+\tilde{\phi}$
is a solution of (\ref{eq:Prob-2}). 
\end{rem*}
To conclude the proof it lasts to find a pair $(\bar{\lambda},\bar{q})$
which is a critical point for $I_{\varepsilon_{1},\varepsilon_{2}}(\lambda,q)$. 

Let us call $G(\lambda,q):=\lambda\alpha(q)B+\lambda^{2}\varphi(q)$.
We have that $\alpha(q)B$ is strictly positive on $\partial M$,
by our assumptions, while by Lemma \ref{lem:phi<0}, $\varphi$ is
strictly negative on $\partial M$. At this point there exists a compact
set $[a,b]\subset\mathbb{R}^{+}$ such that the function $G$ admits
an absolute maximum in $(a,b)\times\partial M$, which also is the
absolute maximum value of $G$ on $\mathbb{R}^{+}\times\partial M$.
This maximum is also $C^{0}$-stable, in the sense that, if $(\lambda_{0},q_{0})$
is the maximum point for $G$, for any function $f\in C^{1}([a,b]\times\partial M)$
with $\|f\|_{C^{0}}$ sufficiently small, then the function $G+f$
on $[a,b]\times\partial M$ admits a maximum point $(\bar{\lambda},\bar{q})$
close to $(\lambda_{0},q_{0})$. By the $C_{0}$ stability of this
maximum $(\lambda_{0},q_{0})$, and by Lemma \ref{lem:espansione},
given $\varepsilon_{1}$ sufficiently small (and $\varepsilon_{2}=o(\varepsilon_{1}^{2})$),
there exists a pair $\left(\lambda_{\varepsilon_{1}},q_{\varepsilon_{1}}\right)$
which is a maximum point for $J_{g}\left(\tilde{W}_{\sqrt{\lambda\varepsilon_{1}},q}+\lambda\varepsilon_{1}\tilde{V}_{\sqrt{\lambda\varepsilon_{1}},q}\right)$,
and, in turn, that there exists a pair $\left(\bar{\lambda}_{\varepsilon_{1}},\bar{q}_{\varepsilon_{1}}\right)$
which is a maximum point for $I_{\varepsilon_{1},\varepsilon_{2}}\left(\lambda,q\right)$.
This implies, in light of the above Remark, that $\tilde{W}_{\sqrt{\bar{\lambda}_{\varepsilon_{1}}\varepsilon_{1}},\bar{q}_{\varepsilon_{1}}}+\bar{\lambda}_{\varepsilon_{1}}\varepsilon_{1}\tilde{V}_{\sqrt{\bar{\lambda}_{\varepsilon_{1}}\varepsilon_{1}},\bar{q}_{\varepsilon_{1}}}+\tilde{\phi}$
is a solution of (\ref{eq:Prob-2}), and the proof for the case $\alpha>0$
is complete.

For the case $\beta>0$ we choose 
\[
\delta=\lambda\varepsilon_{2}^{\frac{1}{3}}\text{ and }\varepsilon_{1}=o(\varepsilon_{2}^{\frac{2}{3}})
\]
in order to have 
\[
J_{g}\left(\tilde{W}_{\lambda\varepsilon_{2}^{\frac{1}{3}},q}+\lambda^{2}\varepsilon_{2}^{\frac{2}{3}}\tilde{V}_{\delta,q}\right)=A+\varepsilon_{2}^{\frac{4}{3}}\left(\lambda\beta(q)C+\lambda^{4}\varphi(q)\right)+o(\varepsilon_{2}^{\frac{4}{3}}),
\]
and the proof follows identically.
\end{proof}
\begin{rem}
\label{rem:sign-changing}We give an example of sign changing perturbation
$\alpha(q)$ such that problem (\ref{eq:Prob-2}) admits blowing up
sequences of solutions. Since $\partial M$ is compact, there exists
a $q_{0}\in\partial M$ maximum point for $\varphi$. We take a $\alpha\in C^{2}(\partial M)$
which has a positive local maximum in $q_{0}$, and that is negative
somewhere. We choose 
\begin{align*}
\delta & =\sqrt{\lambda\varepsilon_{1}}\\
\varepsilon_{2} & =o(\varepsilon_{1}^{2})
\end{align*}
as in the previous proof. By construction, the pair $(\lambda_{0},q_{0})=\left(-\frac{B\alpha(q_{0})}{2\varphi(q_{0})},q_{0}\right)$
is a $C^{0}$-stable critical point for $G(\lambda,q)$, infact $\nabla_{\lambda,q}G(\lambda_{0},q_{0})=0$
and the Hessian matrix is negative definite. Then we can repeat the
arguments of Theorem \ref{thm:main2}. The costruction of a sign changing
$\beta$ is completely analogous. 
\end{rem}

\end{document}